\documentclass[a4paper,11pt,epsfig,graphicx,subfigure,pstricks,colors,twoside]{article}
\usepackage[dvips]{graphics,pstricks}
\usepackage[english]{babel}
\usepackage{epsf}
\usepackage{verbatim}
\usepackage{stmaryrd}
\usepackage{latexsym}
\usepackage{makeidx}
\usepackage{epsfig}
\usepackage{amsfonts}
\usepackage{amsmath}
\usepackage{amssymb}

\input xy
\xyoption{all}
\makeindex
\large \normalsize
\newtheorem{Def}{Definition}[section]
\newtheorem{Prop}[Def]{Proposition}
\newtheorem{Teo}[Def]{Theorem}

\newtheorem{Lem}[Def]{Lemma}
\def\norm#1{\left |#1\right |}

\def\I{\tilde I}

\def\epsilon{\varepsilon}
\def\endproof{\ \hfill\hbox{\vbox{\hrule\hbox{\vrule
height5pt\kern5pt\vrule height5pt}\hrule}}\par\medskip\rm}

\title{\bf 
A new estimate on Evans' Weak KAM approach
}
\date{May 2, 2012}
\author{ \sc Olga Bernardi \qquad Franco Cardin \qquad Massimiliano Guzzo \\ \\
Universit\`a degli Studi di Padova\\
Dipartimento di Matematica \\ Via Trieste, 63 - 35121 Padova, Italy \\ \\
obern@math.unipd.it \ cardin@math.unipd.it \ guzzo@math.unipd.it}
\begin{document}
\maketitle
\begin{abstract} \par \noindent We consider a recent approximate variational principle for weak KAM theory proposed by Evans. 
As in the case of classical integrability, 
for one dimensional mechanical Hamiltonian systems all the computations can be carried out 
explicitly. In this setting, we illustrate the geometric content of the theory and prove new lower bounds 
for the estimates related to its dynamic interpretation. These estimates also extend to the case of $n$ degrees of freedom. 
\end{abstract}

\section{Introduction} The integrability of classical
mechanical systems follows from the existence of regular global solutions to 
the Hamilton-Jacobi equation
\begin{equation} \label{HJintro} 
H(\tilde p + {\partial u\over \partial q}, q) = \bar{H}(\tilde p) 
\end{equation}
where both the generating function $u(\tilde p,q)$ and the 
Hamiltonian $\bar{H}(\tilde p)$  are unknown. 
It is well known (by the Liouville-Arnol'd Theorem) that global solutions to this problem exist only for a very special class of 
mechanical systems, namely, those having a complete set of first integrals. Although most mechanical systems are not integrable 
in this sense, many are quasi-integrable, that is they have the form
\begin{equation}
H(I,\varphi )=h(I)+ \varepsilon f(I,\varphi ),
\label{ham1}
\end{equation}
where $(I,\varphi) \in {\Bbb R}^n \times {\Bbb T}^n$ are action-angle
variables. The new approach to Hamiltonian perturbation theories motivated by Poincar\'e contributions
culminated with the celebrated KAM theorem  \cite{Arno63}, \cite{kolmo54}, \cite{Moser58}. \\

\indent
Since the early 1980's alternative approaches to the study of
non-integrable Hamiltonians based on variational methods and
PDE techniques \cite{Mather1,Mather3,Mather2}, \cite{Mane}, \cite{fathi1,fathi2} have led to the
formulation of the so-called weak KAM theory. Within the Tonelli setting, that is, assuming positive definite superlinear 
Lagrangians and Hamiltonians, the main
results of this theory are the existence of invariant
(action-minimizing) sets, generalizing KAM tori, and the existence of
global weak solutions to the Hamilton-Jacobi equation (\ref{HJintro}). In particular, it has been proved 
in various contexts (homogenization \cite{LPV}, variational and viscosity \cite{fathibook})
that if $H(I,\varphi)$ is Tonelli, then for 
any $I$ the Hamilton-Jacobi problem
\begin{equation} \label{HJintroIp} 
H(I + {\partial u\over \partial \varphi}, \varphi) = \bar{H}(I) 
\end{equation}
admits Lipschitz continuous solutions, with ``effective'' Hamiltonian
\begin{equation} \label{eff}
\bar{H}(I) := \inf_{u \in C^1({\Bbb T}^n)} \max_{\varphi \in {\Bbb T}^n}
H(I + {\partial u\over \partial \varphi}, \varphi).\end{equation}
In the terminology of \cite{fathi1}, these solutions are called weak KAM. Most of the 
dynamic interpretations of these weak solutions have been related to 
Aubry-Mather theory (see for example \cite{contritu}, \cite{kalo} and the
references therein). \\ 
\indent The starting point of the  present paper is an
innovative formulation of weak KAM theory given by Evans \cite{E1,E2}. The main 
outcome of this new variational construction, inspired by Aronsson's variational principle, is a sequence of smooth functions $u_k(\tilde I,\varphi)$ which define, 
for any value of the parameters ${\tilde I}$, a dynamics and a density measure $\sigma_k({\tilde I},\varphi)$ 
on the torus ${\Bbb T}^n$. The convergence of this torus dynamics to a linear flow is expressed precisely through 
the asymptotic formula (2.22) in \cite{E2}. Moreover, an estimate of 
how the torus flow approximates the genuine Hamiltonian flow of $H(I,\varphi)$ 
is expressed through the asymptotic formula (2.21) of 
\cite{E2}.  We refer to Section \ref{SEZIONE 2} for complete details. 
The fundamental relations (2.21) and (2.22) of \cite{E2} are expressed in the 
form of upper bounds. \\
\indent The first goal of this paper is to offer a detailed geometric and dynamic representation, summarized in Figures \ref{fig1} and 
\ref{fig2}, of several evolutions and flows arising from Evans framework. Moreover, we complete the fundamental estimates of Evans, (2.21) and (2.22),  by also measuring  the gap $d_k$ 
between the original Hamiltonian flow and a crucial approximate dynamics introduced by Evans. \\
\indent We remark that in the generic $n$ dimensional case, there 
exists no  explicit expression for the $u_k(\tilde{I},\varphi)$, although
numerical approximations may be obtained via a finite difference 
scheme \cite{falconer}. There is one case in which the sequences 
$u_k$ have an explicit analytic representation and a simple mechanical
interpretation, namely, the case of one degree of freedom. In particular, in this paper we show that for a mechanical system
\begin{equation} 
H(I,\varphi) := {I^2 \over 2}+f(\varphi ) ,
\label{mechamn1}
\end{equation}
with $(I,\varphi) \in {\Bbb R}^+\times \mathbb{S}^1$, 
the functions $u_k$ are precisely the solutions of the   
Hamilton-Jacobi equations for the modified Hamiltonians
\begin{equation} \label{modHam}
H_k(I,\varphi) =  \frac{I^2}{2} + f(\varphi) + \frac{1}{k} \log I.
\end{equation}
As a consequence of the special form of the term $\frac{1}{k}
\log I$, for any $I\in {\Bbb R}^+$ the solution $u_k$ of the
Hamilton-Jacobi equation is explicit up to quadratures of elementary functions 
and the special Lambert function. By taking advantage of this 
explicit analytic expression for the $u_k$, we can prove better convergence 
properties than the more general ones given in \cite{E2}, give new lower bounds in the  
inequalities (2.21), (2.22) of  \cite{E2} and also exhibit an explicit example of  
singular convergence of the measures $\sigma_k$. \\
These new lower bounds constrain  the $\sigma_k$-convergence of the approximate dynamics to a linear flow 
to be, in general,  no faster than $1/k^2$. In Section \ref{quattro}, we see that these one dimensional estimates also have further consequences 
in the integrable $n$ dimensional case, with $n-1$ ignorable variables. \\ 
\indent We also remark that the present one dimensional study may provide a basis 
for a perturbation approach to single resonances in Hamiltonian systems, whose 
normal forms are represented by perturbations of the mechanical pendulum. \\
\indent The paper is organized as follows. In Section \ref{SEZIONE 2} we review the fundamentals 
of the Evans theory and we offer a geometric and dynamic representation of several evolutions and flows of this framework. Moreover, we measure the gap between the original Hamiltonian flow and a crucial approximate dynamics introduced by Evans. 
Section \ref{sezione 3}  is devoted to explicit solutions and convergence results in the one dimensional case. In Section \ref{quattro},
 by exploiting our explicit knowledge of the sequences $u_k$ and $\sigma_k$ in the one dimensional case, we first propose refined 
asymptotic estimates for the integrals involved in formulas  
(2.21), (2.22) of \cite{E2} --integrals (\ref{estsigmakeff}) and (\ref{estsigmak}) here-- and then we discuss some consequences 
in the quadrature-integrable $n$ dimensional case. 
Sections \ref{cinque}, \ref{sei} and \ref{sette} are devoted to the proofs. 
In Section \ref{W} we 
review some properties of the special Lambert function. 
\section{Dynamic picture of Evans theory} \label{SEZIONE 2}
In \cite{E1,E2} Evans introduces a new variational version of weak KAM 
theory, whose outcome is a sequence of functions $u_k({\tilde I},\varphi)$  which define, for any value of the 
parameter ${\tilde I} \in \mathbb{R}^n$ and any index $k\in {\Bbb N}$, a
dynamics and a density measure $\sigma_k({\tilde I},\varphi)$ on the 
torus ${\Bbb T}^n$. The properties of this torus dynamics and its relations 
with the original Hamiltonian 
flow represent the dynamic interest of the theory. \\
\indent More precisely, instead of looking for minimizers $u(\tilde{I},\varphi)$ for the 
sup-norm of $H(\tilde{I} + {\partial u\over \partial \varphi}, \varphi)$ over 
$\mathbb{T}^n$, Evans looks for minimizers $u_k(\tilde{I},\varphi)$ of the functional
\begin{equation} \label{functional}
I_k[u] := \int_{\mathbb{T}^n} e^{kH(\tilde{I}+{\partial u\over \partial \varphi}
  ,\varphi)} d\varphi.
\end{equation}
Under suitable hypotheses\footnote{
Precisely, $H$ is periodic in the $\varphi$ variables; 
$H$ is convex in the $I$ variables; there exists $C > 0$ such that, for 
any $I\in {\Bbb R}^n$ and $\varphi\in {\Bbb T}^n$:
$\max \Big \{ \norm{{\partial^2 H\over 
\partial I^2}},  {\norm{{\partial^2 H\over 
\partial I\partial \varphi}}\over 1+\norm{I}},{\norm{{\partial^2 H\over 
\partial \varphi^2}}\over 1+\norm{I}^2}\Big \}\leq C$.}
on $H$, the minimizers $u_k$   
turn out to be smooth and uniquely defined when one  requires that
$\int_{\mathbb{T}^n} u_k d\varphi = 0$. After defining the density measure 
over ${\Bbb T}^n$
\begin{equation} \label{sigma}
\sigma_k(\tilde{I},\varphi) := e^{k(H(\tilde{I} + {\partial u_k \over \partial \varphi},\varphi) 
- \bar{H}_k(\tilde{I}))}
\end{equation}
where
\begin{equation} \label{acca}
\bar{H}_k(\tilde{I}) := \frac{1}{k} \log  \int_{\mathbb{T}^n} e^{k H(\tilde{I} +
{\partial u_k\over \partial \varphi},\varphi)} d\varphi  ,
\end{equation}
Evans  (Theorems 2.1 and 3.1 in \cite{E1}) proves that
\begin{equation}\label{intestimate}
\lim_{k \to +\infty} \int_{\mathbb{T}^n} H(\tilde{I} + {\partial u_k
\over \partial \varphi},\varphi) \sigma_k(\tilde{I},\varphi) d\varphi = 
\bar{H}(\tilde{I})=\lim_{k \to +\infty} \bar{H}_k(\tilde{I}) ,
\end{equation}
where $\bar{H}$ is the usual effective Hamiltonian --see (\ref{eff})-- of weak KAM theory. \\
\indent Since the functions $u_k$ are smooth, they may be used to generate canonical
transformations $(I,\varphi)\mapsto (\tilde I,\tilde \varphi)$
up to the inversion of\footnote{Of course, the functions 
$U_I,U_\varphi$ depend on the parameter $k$. But, since we will not    
define their limit for $k$ tending to infinity, we prefer to use this
simplified implicit notation.}
\begin{eqnarray} \label{genfun1}
I = U_I(\tilde I,\varphi), \quad U_I(\tilde I,\varphi)=
\tilde I+{\partial u_k\over \partial \varphi}(\tilde I,\varphi)
\end{eqnarray}
and
\begin{eqnarray} \label{genfun2}
\tilde \varphi=U_\varphi (\tilde I,\varphi), \quad  
U_\varphi (\tilde I,\varphi)=\varphi+{\partial u_k\over \partial I}
(\tilde I,\varphi).
\end{eqnarray} 
For every fixed  $\tilde I\in {\Bbb R}^n$, Evans introduces the dynamics on 
the torus $\mathbb{T}^n$ by the differential equation
\begin{equation}
\dot \varphi = {\partial H\over \partial I}\Big (
\tilde I+{\partial u_k\over \partial \varphi}(\tilde I,\varphi),\varphi \Big),
\label{torusdyn}
\end{equation}
whose flow will be here denoted  by ${\cal C}^t_{\tilde I}(\varphi)$. 
This torus flow ${\cal C}^t_{\tilde  I}(\varphi)$ preserves the measure 
defined by $\sigma_k(\tilde{I},\varphi)$, see \cite{E2}. Indeed, from the Euler-Lagrange 
equation related to the variation of $I_k[u]$,  one obtains
\begin{equation}
{\rm Div}_\varphi \Big (\sigma_k(\tilde{I},\varphi) 
 {\partial H\over \partial I}\Big (
\tilde I+{\partial u_k\over \partial \varphi}(\tilde I,\varphi),\varphi \Big )
\Big ) =0 .
\label{divsigma}
\end{equation}
Starting from  $\mathcal{C}_{\tilde{I}}^t(\varphi)$, 
and inspired by equations (2.16)-(2.18) of \cite{E2}, 
we define the following two evolutions
\begin{equation} \label{phi}
\Phi^t: (\tilde{I},\varphi) \to (I^t,\varphi^t) := (U_I(\tilde{I},\mathcal{C}_{\tilde{I}}^t(\varphi)),\mathcal{C}_{\tilde{I}}^t(\varphi))
\end{equation}
and
\begin{equation} \label{phi tilde}
\tilde{\Phi}^t: (\tilde{I},\varphi) \to (\tilde{I}^t,\tilde{\varphi}^t) := (\tilde{I},U_{\varphi}(\tilde{I},\mathcal{C}_{\tilde{I}}^t(\varphi))),
\end{equation}
which are obtained as the composition of the flow 
$$\tau^t: (\tilde{I},\varphi) \to (\tilde{I}, {\cal C}^t_{\tilde I}(\varphi))$$ 
and the transformations (\ref{genfun1}), (\ref{genfun2}) of actions and angles respectively
(see Figure \ref{fig1}). We remark that the  
$(I^t,\varphi^t)$ and  $(\tilde{I}^t,\tilde{\varphi}^t)$ 
are not necessarily conjugate, since
expressions (\ref{genfun1}), (\ref{genfun2}) are not necessarily invertible. 
\begin{figure}[!t]
$$\xymatrix{& (I^t, \varphi^t) \\
(\tilde{I},\varphi) \ar[r]^{\tau} \ar[dr]_{\tilde{\Phi}^t} \ar[ur]^{\Phi^t} &
  (\tilde{I},\mathcal{C}^t_{\tilde{I}}(\varphi)) \ar[d]^{\textit{angles lift}}
  \ar[u]_{\textit{actions lift}}\  \\ 
& (\tilde{I}^t,\tilde{\varphi}^t)}$$
\caption{\small{The dynamics (\ref{phi}) and (\ref{phi tilde}) represent different 
lifts of the torus flow dynamics. We note that
in both cases  the domain  is the mixed variables set ${\Bbb R}^n\times 
{\Bbb T}^n$.}}
\label{fig1}
\end{figure}
\begin{figure} 
$$\xymatrix{
(I^t_H,\varphi^t_H) \ar@{.}[r]^{(\star)} & (I^t,\varphi^t) \\
(I,\varphi) \ar[u]^{\Phi^t_H} & 
\hbox{$(\tilde{I},\varphi)$} \ar@{-->}[l]^{U_I \times
      \textit{ id}} \ar[u]_{\Phi^t} \ar@{-->}[r]^{\textit{id } \times U_{\varphi}}
    \ar[d]_{\tilde{\Phi}^t} & (\tilde{I},\tilde{\varphi})
    \ar[d]^{\Phi^t_{\bar{H}_k}} \\
& (\tilde{I}^t,\tilde{\varphi}^t)
    \ar@{.}[r]_{(\bullet)} & 
(\tilde I^t_{\bar{H}_k},\tilde \varphi^t_{\bar{H}_k})
}$$\caption{\small{The dynamics 
are represented by solid arrows; the transformations of variables by dashed 
arrows; dotted lines link the evolutions (\ref{phi}), 
(\ref{phi tilde}) and the Hamiltonian flows (\ref{flowham}).}}
\label{fig2}
\end{figure}
\subsection{Relations between the torus dynamics and the Hamiltonian flows} \label{relations}
We stress that the very dynamic relevance of the 
$\tilde{I}$-collection of flows
 $\mathcal{C}^t_{\tilde{I}}(\varphi)$ lies in the relation between the
 orbits of  (\ref{phi}), (\ref{phi tilde}) and those of the 
Hamiltonian flows
\begin{equation}
\Phi_H^t: (I,\varphi) \to (I^t_H,\varphi^t_H) \qquad \Phi_{\bar{H}_k}^t: 
(\tilde I,\tilde \varphi) \to (\tilde I^t_{\bar{H}_k},\tilde
\varphi^t_{\bar{H}_k}) := (\tilde I,\tilde \varphi + t 
\frac{\partial \bar{H}_k}{\partial I}(\tilde I))
\label{flowham}
\end{equation}
of $H$ and $\bar{H}_k$ (defined in ({\ref{acca})) respectively. \\
\indent In \cite{E2} the relation between $\tilde \Phi^t$ and the Hamiltonian
flow $\Phi^t_{\bar{H}_k}$ ($(\bullet)$ in Figure \ref{fig2}) 
is expressed through the asymptotic formula (2.22).  Specifically, Evans
proves that there exists a constant $C_R > 0$ such that
\begin{equation}
E_1(k) := \int_{\norm{\tilde I}\leq R} \int_{{\Bbb T}^n}
\Big \vert {\partial {\tilde \varphi^t} \over \partial t} -
{\partial {\bar H}_k\over \partial I} 
(\tilde I){\Big \vert}^2 \sigma_k(\tilde{I},\varphi) d\varphi d\tilde I \leq {C_R\over k}
\label{estsigmakeff}
\end{equation}
$\forall t \in \mathbb{R}$. Moreover he shows (formula (2.21) in \cite{E2}) that for some $C^R > 0$
\begin{equation}
E_2(k) := \int_{{\Bbb T}^n} \Big \vert {\partial  I^t\over \partial t}
 +{\partial H\over \partial\varphi} ( \Phi^t(\tilde I,\varphi) ){\Big \vert}^2 \sigma_k(\tilde{I},\varphi) 
d\varphi \leq {C^R\over k}
\label{estsigmak}
\end{equation}
$\forall t \in \mathbb{R}$ and $|\tilde{I}| \le R$. \\
\indent By the next proposition, we complete the dynamic picture by studying the relation $(\star)$ in Figure \ref{fig2}. 
More precisely, using the estimate  (\ref{estsigmak}) of Evans,   we measure  the gap 
\begin{equation} \label{distanza}
d_k(t,\tilde I,\varphi) := |\Phi^t(\tilde I,\varphi)-\Phi^t_H(U_I(\tilde
I,\varphi),\varphi)| = |(I^t-I^t_H,\varphi^t - \varphi^t_H)|
\end{equation}
between the curves $(I^t_H,\varphi^t_H)$ and $(I^t,\varphi^t)$     in terms of $\sigma_k$.
\begin{Prop}\label{prop}
Let $\lambda_H > 0$ be a Lipschitz constant for the Hamiltonian vector 
field of $X_H$. Then, we have
\begin{equation}
\int_{{\Bbb T}^n} d_k(t,\tilde I,\varphi) 
\sigma_k(\tilde I,\varphi)d\varphi \leq \frac{(e^{\lambda_H t} - 1)(1+C^R)}{
\lambda_H \sqrt{k}}
\label{estd}
\end{equation}
$\forall t \in \mathbb{R}$, $k\in {\Bbb N}$ and $|\tilde{I}| \le R$. In particular,
$$
\lim_{k\rightarrow +\infty}
\int_{{\Bbb T}^n} d_k(t,\tilde I,\varphi)\sigma_k(\tilde I,\varphi)d\varphi
=0
$$
$\forall t \in \mathbb{R}$ and $|\tilde{I}| \le R$.
\end{Prop}
We remark that the presence of the exponential 
term $e^{\lambda_H t}$ in (\ref{estd}) is not 
surprising, since any small correction to a differential equation typically produces 
an exponential divergence of the solutions. In Section \ref{sette} we provide some further
detail regarding this divergence, by considering an example with exponential 
divergence due to the presence of a hyperbolic equilibrium point of 
the Hamiltonian flow. 
\vskip 0.4 cm
\noindent The proof of the proposition is based on the next technical 
\begin{Lem}\label{lemgw}
For any $t \in \mathbb{R}$, $k\in {\Bbb N}$, $\tilde I\in {\Bbb R}^n$ and $\varphi\in {\Bbb T}^n$, we have
\begin{equation}
d_k(t,\tilde I,\varphi) \leq  e^{\lambda_H t}
\int_0^t \norm{{\partial I^t\over \partial t}\Big{|}_{t = s} +{\partial H\over \partial \varphi}(\Phi^s(\tilde I,\varphi)) }e^{-\lambda_H s}ds.
\label{ineqdk}
\end{equation}
\end{Lem}
\noindent
\textit{Proof of Lemma.}
The time derivatives 
of $I^t,\varphi^t$ satisfy (see (\ref{phi}))
\begin{eqnarray*} 
{\partial I^t\over \partial t} &=& \frac{\partial^2 u_k}
{\partial \varphi^2}(\tilde{I},\varphi^t)\frac{\partial H}{\partial I}
(U_I(\tilde I,\varphi^t),\varphi^t)=
\frac{\partial^2 u_k}
{\partial \varphi^2}(\tilde{I},\varphi^t)\frac{\partial H}{\partial I}
(I^t,\varphi^t),\cr\cr
{\partial \varphi^t\over \partial t} &=& 
\frac{\partial H}{\partial I}(U_I(\tilde I,\varphi^t),\varphi^t)=
\frac{\partial H}{\partial I}(I^t,\varphi^t)
\end{eqnarray*}
so that the functions $I^t,\varphi^t$ may be interpreted as the solutions 
of the following $\tilde I$-parametric differential equation
\begin{equation*} 
\begin{cases}
\dot I = \frac{\partial^2 u_k}
{\partial \varphi^2}(\tilde{I},\varphi)\frac{\partial H}{\partial I}
(I,\varphi) \\
\dot \varphi = 
\frac{\partial H}{\partial I}(I,\varphi)
\end{cases}
\end{equation*}
with special initial conditions $(I,\varphi) = (U_I(\tilde I,\varphi),\varphi)$, and solutions denoted by  $\Phi^t(\tilde{I},\varphi)$, see (\ref{phi}). \\
Let
$$X_{\tilde I}(I,\varphi) = 
\Big( \frac{\partial^2u_k}{\partial
  \varphi^2}(\tilde{I},\varphi)\frac{\partial H}{\partial I}(I,\varphi)
,\frac{\partial H}{\partial I} (I,\varphi)\Big)$$
and
$$X_H(I,\varphi) = 
\Big( -\frac{\partial H}{\partial \varphi}(I,\varphi),
\frac{\partial H}{\partial I} (I,\varphi)\Big)  .$$
Since 
$$
d_k(t,\tilde{I},\varphi) = | (I^t-I^t_H,\varphi^t-\varphi^t_H ) |,
$$
its time derivative 
\begin{equation*}
{\partial d_k \over \partial t}
= {1\over d_k} \Big (I^t-I^t_H,\varphi^t-\varphi^t_H\Big )\cdot 
\Big (X_{\tilde I}(I^t,\varphi^t)- X_H(I^t_H,\varphi^t_H)\Big )
\end{equation*}
is well defined only for $d_k(t,\tilde I,\varphi)>0$, for example, $d_k(0,\tilde{I},\varphi) = 0$. 
In order to overcome the lack of differentiability, for a constant $\epsilon>0$ we consider the function
$$
x_k(t,\tilde I,\varphi ) := \sqrt{d_k(t,\tilde I,\varphi )^2+\epsilon^2},
$$
whose time derivative
\begin{equation*}
{\partial x_k \over \partial t}
= {1\over x_k} \Big (I^t-I^t_H,\varphi^t-\varphi^t_H\Big )\cdot 
\Big (X_{\tilde I}(I^t,\varphi^t)- X_H(I^t_H,\varphi^t_H)\Big )
\end{equation*}
satisfies
\begin{eqnarray*}
&\norm{{\partial x_k \over \partial t}} \le
{d_k \over x_k} \norm{X_{\tilde I}(I^t,\varphi^t)- X_H(I^t_H,\varphi^t_H)}\leq 
\norm{X_{\tilde I}(I^t,\varphi^t)- X_H(I^t_H,\varphi^t_H)} \le &\cr
&\leq \norm{X_{\tilde I}(I^t,\varphi^t)-X_H(I^t,\varphi^t)} +
\norm{X_H(I^t,\varphi^t)- X_H(I^t_H,\varphi^t_H)}.
\end{eqnarray*}
As a consequence, we have ($d_k < x_k$)
$$
\norm{{\partial x_k \over \partial t}} \le \norm{X_{\tilde I}(I^t,\varphi^t)-X_H(I^t,\varphi^t)}+\lambda_H d_k <
\norm{X_{\tilde I}(I^t,\varphi^t)-X_H(I^t,\varphi^t)}+\lambda_H x_k,
$$
where $\lambda_H > 0$ is a Lipschitz constant for $X_H$. By recalling now the classical \textit{a priori} estimate\footnote{In brief, the \textit{a priori} upper bound estimate lemma, see \cite{SC}: if $|f'(t)| < M(t, |f(t)|)$ 
and $g(t)$ solves $\dot{g}(t) = M(t,g(t)) \text{ with } g(0) = f(0)$, then $|f(t)| \le g(t)$.}, we obtain
$$
x_k (t,\tilde I,\varphi) \leq  \int_0^t \norm{X_{\tilde I}(I^s,\varphi^s)-X_H(I^s,\varphi^s)}e^{\lambda_H (t-s)}ds + e^{\lambda_H t} x_k(0),
$$
with $x_k(0)=\epsilon$, so that
$$
d_k (t,\tilde I,\varphi) \leq  e^{\lambda_H t}\Big (
\int_0^t \norm{X_{\tilde I}(I^s,\varphi^s)-X_H(I^s,\varphi^s)}e^{-\lambda_H s}ds
+\epsilon \Big ) .
$$
As a consequence, by considering arbitrarily small $\epsilon > 0$, we conclude that
\begin{equation*}
d_k (t,\tilde I,\varphi) \leq  e^{\lambda_H t}
\int_0^t \norm{X_{\tilde I}(I^s,\varphi^s)-X_H(I^s,\varphi^s)}
e^{-\lambda_H s}ds .
\label{ineqdk0}
\end{equation*}
Finally, since
$$
X_{\tilde I}(I^s,\varphi^s)-X_H(I^s,\varphi^s)=
\Big ( {\partial I^t\over \partial t}\Big{|}_{t=s} +
{\partial H\over \partial \varphi}(I^s,\varphi^s),0\Big )
=\Big ( {\partial I^t\over \partial t}\Big{|}_{t=s}+{\partial H\over \partial \varphi}
(\Phi^s(\tilde I,\varphi)),0\Big ),
$$
we have
\begin{equation*}
d_k (t,\tilde I,\varphi) \leq  e^{\lambda_H t}
\int_0^t \norm{{\partial I^t\over \partial t}\Big{|}_{t=s}+
{\partial H\over \partial \varphi}(\Phi^s(\tilde I,\varphi))}
e^{-\lambda_H s}ds
\end{equation*}
for any $t \in \mathbb{R}$, $k\in {\Bbb N}$, $\tilde I\in {\Bbb R}^n$ and $\varphi\in {\Bbb T}^n$. \hfill $\Box$
\vskip 0.4 cm
\noindent {\it Proof of Proposition \ref{prop}.} By considering the $\sigma_k$-average of $d_k$ and the inequality (\ref{ineqdk}), we obtain
\begin{equation*}
\int_{{\Bbb T}^n} d_k(t,\tilde{I},\varphi) \sigma_k(\tilde{I},\varphi) d\varphi \leq  e^{\lambda_H t}
\int_0^t e^{-\lambda_H s} \int_{{\Bbb T}^n}
\norm{{\partial I^t\over \partial t}\Big{|}_{t=s}+{\partial H\over \partial
    \varphi}(\Phi^s(\tilde I,\varphi)) }\sigma_k(\tilde I,\varphi)d\varphi ds.
\label{ineqdk2}
\end{equation*}
Since now for any $x\geq 0$ and $k \in \mathbb{N}$
$$
x \leq \max \Big ( {1\over \sqrt{k}}, \sqrt{k}\ x^2 \Big ) \leq 
{1\over \sqrt{k}}+\sqrt{k}\  x^2,
$$
we conclude that
\begin{eqnarray*}
\int_{{\Bbb T}^n}
\norm{{\partial I^t\over \partial t}\Big{|}_{t=s}+{\partial H\over \partial
    \varphi}(\Phi^s(\tilde I,\varphi)) }\sigma_k(\tilde{I},\varphi) d\varphi
& \leq & \int_{{\Bbb T}^n} 
\Big ({1\over \sqrt{k}}+
\sqrt{k} \norm{{\partial I^t\over \partial t}\Big{|}_{t=s} + {\partial H\over \partial
    \varphi}(\Phi^s(\tilde I,\varphi)) }^2\Big )
\sigma_k(\tilde{I},\varphi) d\varphi \cr
& \leq &
{1\over \sqrt{k}}+{C^R\over \sqrt{k}},
\end{eqnarray*}
where, in the last inequality, we have used the estimate (\ref{estsigmak}) of Evans. Since the right hand side of the inequality does not depend on $s$, 
we immediately obtain (\ref{estd}). \hfill $\Box$
\section{Explicit solutions and convergences in the one dimensional case} \label{sezione 3}
This section presents the explicit formula for the minimizers $u_k$ of the functional $I_k[u]$ defined in (\ref{functional}) for the Hamiltonian systems (\ref{mechamn1}) with  one degree of freedom. \\ 
\indent  For simplicity we consider only the 
action interval $\tilde I>0$: the case $\tilde I<0$ can be obtained by 
symmetry (see Remark \textbf{(IV)} below). In the 
sequel we make extensive use of the Lambert function $W$, defined implicitly by 
$z=W(z)e^{W(z)}$, and also its asymptotic properties. (We refer the reader to the technical Section \ref{W} and to \cite{cor1},\cite{cor2}). 
\begin{Def} \label{def fond}
For $H(I,\varphi) = I^2/2 + f(\varphi)\in {\cal C}^2({\Bbb R}\times {\Bbb
  S}^1)$, let us define: 
\begin{itemize}
\item[$(i)$] For $\I >0$, the sequence of functions $c_k(\I )\in {\Bbb R}$ by inversion of
\begin{equation}
c \longmapsto \I= \frac{1}{2 \pi} \int_0^{2\pi} \gamma_k(c,\varphi )d\varphi  ,
\label{ark}
\end{equation}
where
\begin{equation} 
\gamma_k(c,\varphi) := \sqrt{\frac{W(e^{2(c-f(\varphi))k}k)}{k}}  .
\label{gkdin}
\end{equation}
\item[$(ii)$] For $\I >0$, the function $c(\I)\in {\Bbb R}$ by inversion of
\begin{equation}
c \longmapsto \I= \frac{1}{2 \pi} \int_{0}^{2\pi} 
\gamma_0(c,\varphi)d\varphi  ,
\label{ark0}
\end{equation}
where
\begin{equation} \label{gammazero}
\gamma_0(c,\varphi) := 
\begin{cases} \sqrt{2(c-f(\varphi))} & \textnormal{if } c > f(\varphi) \\
0 & \textnormal{otherwise}
\end{cases}
\end{equation}
\item[$(iii)$] For $(\I,\varphi)\in (0,+\infty) \times {\Bbb S}^1$, the sequences of functions
\begin{equation} \label{soluzione}
u_k(\I,\varphi) := \I(\pi-\varphi)- \frac{1}{2\pi} \int^{2\pi}_0 
\int^y_0 \gamma_k(c_k(\I),x)
dx dy + \int^{\varphi}_0 \gamma_k(c_k(\I),x)dx 
\end{equation}
and
\begin{equation} \label{sigmak}
\sigma_k(\I,\varphi) := {{1\over \gamma_k(c_k(\I),\varphi)}\over 
\int_0^{2\pi}{1\over \gamma_k(c_k(\I),x)}dx} .
\end{equation}
\item[$(iv)$] For $(\I,\varphi)\in (0,+\infty) \times {\Bbb S}^1$, the 
function
\begin{equation} \label{soluzione0}
u_0(\I,\varphi) := \I(\pi-\varphi)- \frac{1}{2\pi} \int^{2\pi}_0 
\int^y_0 \gamma_0(c(\I),x)
dx dy + \int^{\varphi}_0 \gamma_0(c(\I),x)dx 
\end{equation}
\item[$(v)$] For $\I>0$, the sequence of functions
\begin{equation} \label{hamiltoniana}
\bar{H}_k(\I) := c_k(\I) + \frac{1}{k} \log \int_{0}^{2\pi} \frac{1}{\gamma_k(c_k(\I),\varphi)}d\varphi  .
\end{equation}
\end{itemize}
\end{Def} 
The convergence properties of the objects defined above are stated in the 
following  
\begin{Teo} \label{prop1} Let us consider $H(I,\varphi) = I^2/2 + 
f(\varphi)\in {\cal C}^2({\Bbb R}\times {\Bbb S}^1)$.  
\begin{itemize}
\item[$(i)$] For any $\I>0$, the functions $u_k(\I,\varphi)$ defined 
in (\ref{soluzione}) are smooth, have zero average and solve the
Euler-Lagrange equation (\ref{divsigma}) for $I_k[u]$. Moreover,
any $u_k(\I,\varphi)$ converges uniformly to $u_0(\I,\varphi)$ on ${\Bbb S}^1$. 
\item[$(ii)$] The functions $\gamma_k(c,\varphi)$ 
defined in (\ref{gkdin}) are smooth and uniformly converging to 
$\gamma_0(c,\varphi)$ defined in (\ref{gammazero}) on $(-\infty,c_*] \times
\mathbb{S}^1$, for any fixed $c_*\in {\Bbb R}$.
\item[$(iii)$] The functions $c_k(\I)$ are pointwise converging to $c(\I)$.
\item[$(iv)$] The functions $\sigma_k(\I,\varphi)$ and 
$\bar{H}_k(\I)$, defined in (\ref{sigmak}) and (\ref{hamiltoniana}) 
respectively, satisfy
\begin{equation}
\lim_{k\rightarrow +\infty}\int_0^{2\pi} 
H (\I+{\partial u_k\over \partial \varphi}, \varphi)
\sigma_k(\I,\varphi)d\varphi= {\bar H}(\I) = \lim_{k \rightarrow +\infty}
\bar{H}_k(\I),
\label{liminthams}
\end{equation}
where
\begin{equation} \label{H-lim}
\bar{H}(\I) := \begin{cases} c(\I) & \textnormal{if } c(\I) > \max f \\ 
\max f & \textnormal{otherwise}
\end{cases}
\end{equation}
\end{itemize}
\end{Teo}
\noindent{\bf Remarks} 
\begin{itemize}
\item[{\bf (I)}] As we will see in Section \ref{cinque uno}, the functions $\gamma_k(c,\varphi)$
 parametrize the level curves for the Hamiltonians (\ref{modHam}) of value $c$ and are well 
defined for any ${\Bbb \varphi}\in {\Bbb S}^1$.  In other words, these level curves
project injectively on $\mathbb{S}^1$. 
Therefore, the action $\I$ in (\ref{ark}) is proportional to the area of the 
phase-space $(0,+\infty)\times {\Bbb S}^1$ under the graph of 
$\gamma_k(c_k(\I),\varphi)$. More precisely, we have
\begin{equation}
\I={1\over 2\pi}\int_0^{2\pi}\gamma_k(c_k(\I),\varphi)d\varphi
\label{actik}
\end{equation}
as well as
\begin{equation}
\I={1\over 2\pi}\int_0^{2\pi}\gamma_0(c(\I),\varphi)d\varphi .
\label{acti0}
\end{equation} 
Let us remark that $\I>0$ corresponds to $c(\I)>\min f$. 
\item [{\bf (II)}]
The functions $u_k(\I,\varphi)$ are smooth 
solutions of the Hamilton-Jacobi equation for (\ref{modHam}),
\begin{equation}\label{HJ} 
\frac{1}{2} (\I + {\partial u_k\over \partial \varphi})^2 + 
f(\varphi) + \frac{1}{k} \log(\I+ {\partial u_k\over \partial \varphi}) = c
\end{equation}
with $c=c_k(\I)$. 
The PDE (\ref{HJ}), at variance with the Hamilton-Jacobi  
equation for (\ref{mechamn1}), admits smooth solutions defined 
over ${\Bbb S}^1$ 
for all values $c>0$. Once again, this follows because  all level curves 
of (\ref{HJ}) project injectively on $\mathbb{S}^1$. \\
Let us also remark that, while in the general $n$ dimensional setting 
Evans (\cite{E1}, Lemma 2.1)  assumes the uniform convergence 
of the sequence $u_k$, passing if 
necessary to a subsequence, in  the one dimensional case we can prove the stronger 
uniform convergence of $u_k$ to $u_0$ on ${\Bbb S}^1$.
\item[{\bf (III)}] For $\I > 0$ such that $c(\I) > \max f$, the function $\varphi \mapsto \I \varphi + u_0(\I,\varphi)$ provides  
a regular solution to the Hamilton-Jacobi equation for the Hamiltonian $H$, see
(\ref{mechamn1}), on the energy level $c(\I)$. Notice that $I \varphi + u_0(\I,\varphi)$ represents the generating function conjugating $H$ to $\bar{H} = c(\I)$. Otherwise, for $c(\I) \le \max f$ the picture differs from the classical integration of one dimensional Hamiltonian systems, because 
$\gamma_0(c(\I),\varphi)$ has angular points for $c(\I)=f(\varphi)$ and the limit 
function $u_0(\I,\varphi)$ is therefore only Lipschitz.
\item[{\bf (IV)}] The case $\I<0$ is obtained via the choice
$$
u_k(\I,\varphi) := \I(\pi-\varphi)+\frac{1}{2\pi} \int^{2\pi}_0 
\int^y_0 \gamma_k(c_k(\vert \I\vert ),x)
dx dy - \int^{\varphi}_0 \gamma_k(c_k(\vert \I\vert ),x)dx  .
$$
As a consequence one also has $\bar{H}_k(\I)=\bar{H}_k(\vert \I\vert)$ 
and $\sigma_k(\I,\varphi)=\sigma_k(\vert \I\vert,\varphi)$.
\end{itemize}
\noindent We devote here some attention to the convergence properties 
of the density measures $\sigma_k$.  In the generic $n$ dimensional case, 
Evans \cite{E1,E2} discusses the consequences of the 
convergence $$\sigma_k \rightharpoonup \sigma \qquad \text{weakly as measures
  on } \mathbb{T}^n$$
possibly through a sub-sequence.

A particularly interesting case corresponds to the convergence of 
the $\sigma_k$ to singular measures on the torus ${\Bbb  T}^n$. Unfortunately, 
the theory of \cite{E1} and \cite{E2} does not provide explicit examples.
In the one dimensional case, if $c(\I) > \max f$,  
the limit of $\sigma_k$ obviously defines a regular measure
on $\mathbb{S}^1$. The case $c(\I) \le \max f$ is actually more
tricky to manage. The following proposition gives  an example of convergence to a singular 
measure:
\begin{Prop} \label{sigma-lim} Let $f(\varphi) = -\cos \varphi$ and 
$\I >0$ be such that $c(\I) = f(\pi)=1$. For any test function $u \in C^{\infty}([0,2\pi);\mathbb{R})$, one has
$$\lim_{k \to +\infty} \int^{2\pi}_0 u(\varphi) \sigma_k(\I,\varphi) d\varphi = u(\pi) .$$
\end{Prop}
\section{Lower bounds and outcomes in the $n$ dimensional case} \label{quattro}
By exploiting our explicit knowledge of the sequences $u_k$ and $\sigma_k$
in the one dimensional case, we first propose to give refined asymptotic estimates 
for the integrals (\ref{estsigmakeff}) and (\ref{estsigmak}).  The estimates are also relevant  for the $n$ dimensional case. 
We start with the following:
\begin{Teo}\label{prop2}
Let us consider $H(I,\varphi) = I^2/2 + f(\varphi)\in {\cal C}^2({\Bbb
  R}\times {\Bbb  S}^1)$, where $f$ is a non constant function.   
\begin{itemize}
\item[$(i)$] For any $r>0$ and $R>0$ satisfying $c(R) > \max f+ r$, 
there exist $K > 0$ and $c_R>0$ such that
\begin{equation} \label{stima1}
E_1(k) \ge \frac{c_R}{k^2}
\end{equation}
$\forall k > K$.
\item[$(ii)$] For any $\I >0$ such that $c(\I) > \max f$, there exists $c_{\I} > 0$ such that
\begin{equation} \label{stima2}
\lim_{k \to +\infty} k^2 E_2(k) = c_{\I} .
\end{equation}
\end{itemize}
In particular, we have
\begin{equation*} 
c_R := 4 \pi \int_{\max f + r}^{c(R)} \frac{- a^2_{\frac{3}{2}}(c) + a_{\frac{5}{2}}(c) a_{\frac{1}{2}}(c)}{a^3_{\frac{1}{2}}(c)} dc \ \ 
\end{equation*}
and
$$c_{\I} := \frac{1}{a_{\frac{1}{2}(c(\I))}} \int^{2\pi}_0 \frac{|f'(\varphi)|^2}{\gamma_0^5(c(\I),\varphi)} d\varphi = \frac{1}{a_{\frac{1}{2}(c(\I))}} \int^{2\pi}_0 \frac{|f'(\varphi)|^2}{[2(c(\I) - f(\varphi))]^{5/2}}d\varphi  ,$$
where
$$ {a}_{\delta}(c) := \int^{2\pi}_0 \frac{1}{\gamma_0^{2\delta}(c,\varphi)} d\varphi = 
\int_0^{2\pi} {1\over [2(c - f(\varphi))]^{\delta}} d\varphi .
$$
\end{Teo}
\vskip 0.4 cm
\noindent Theorem \ref{prop2} provides lower bounds which are also significant  
for the generic $n$ dimensional case. Indeed, for the $n$ degrees of freedom 
mechanical Hamiltonians
\begin{equation} \label{mechanical}
H(I_1,\ldots ,I_n,\varphi_1,\ldots ,\varphi_n)= \sum_{j=1}^n {I_j^2\over 2} + f(\varphi_1,\ldots ,\varphi_n) ,
\end{equation}
let us consider Evans construction of the sequences $u_k^{(n)}$, 
$\bar H_k^{(n)}$ and $\sigma_k^{(n)}$, 
as well as the integrals $E_1(k),E_2(k)$. A relevant question 
is that regarding optimality of the upper bounds (\ref{estsigmakeff}) and
(\ref{estsigmak}) proved in Evans paper. Of course, in the
trivial integrable case $H(I) = I^2/2$, both $E_1(k)$ and $E_2(k)$ are zero. However, already 
in the quadrature-integrable case, for example
$$
f (\varphi_1,\ldots ,\varphi_n):= f(\varphi_1),
$$
we have
\begin{equation} \label{quadra-inte}
{c_R\over k^2} \leq E_1(k) \le {C_R\over k} \qquad \text{and} \qquad {c_{\I}\over k^2} \leq E_2(k) \le {C^R\over k} 
\end{equation}
for any $c(R) >\max f$ (see Theorem 4.1) and sufficiently large $k$. In fact 
all the sequences can be constructed by
referring to the one dimensional case
$$u_k^{(n)}(\I_1,\ldots,\I_n,\varphi_1,\ldots,\varphi_n) := 
u_k^{(1)}(\I_1,\varphi_1), \qquad \sigma_k^{(n)}(\I_1,\ldots,\I_n,\varphi_1,\ldots,\varphi_n) := \sigma_k^{(1)}(\I_1,\varphi_1)$$
and
$$\bar{H}_k^{(n)}(\I_1,\ldots,\I_n) = \bar{H}_k^{(1)}(\I_1) + \sum_{j=2}^n 
{\I_j^2\over 2} .
$$
Moreover,  with regard to the mechanical Hamiltonian systems (\ref{mechanical}), the integrals in (\ref{estsigmakeff}) and (\ref{estsigmak}) can be written precisely in the following form (see Theorem \ref{prop2}):
\begin{equation*}
E_1(k) =  
\int_{|\tilde{I}| \le R} \int_{\mathbb{T}^n} \Big \vert
\frac{\partial}{\partial I} \Big[\frac{1}{2}\left(\tilde{I} + \frac{\partial
    u_k^{(n)}}{\partial \varphi}\right)^2 - \bar{H}_k^{(n)}(\tilde{I})\Big] \Big\vert^2 \sigma_k^{(n)}
 d\varphi d\tilde{I}
\end{equation*}
and
\begin{equation*}
E_2(k) = 
\int_{\mathbb{T}^n} \Big \vert \frac{\partial}{\partial \varphi} \Big[
\frac{1}{2}\left(\tilde{I} + \frac{\partial u_k^{(n)}}{\partial \varphi}\right)^2 +
f(\varphi)\Big] \Big\vert^2  \sigma_k^{(n)}d\varphi.
\end{equation*}
As a consequence, we immediately conclude that the estimates (\ref{quadra-inte}) also
apply to the quadrature-integrable $n$ dimensional case.

\section{Proof of Theorem \ref{prop1}} \label{cinque}

\subsection{Explicit formulas for $u_k$} \label{cinque uno}
According to \cite{E2}, the function $u_k$ is a minimizer of the 
functional $I_k[u]$ defined in (\ref{functional}), whose  
Euler-Lagrange equation is
$$\sum_{j=1}^n {\partial \ \over \partial \varphi_j}
\left(e^{kH(\I+{\partial u_k\over \partial \varphi},\varphi)} 
\frac{\partial H}{\partial  I}(\I+{\partial u_k\over \partial \varphi},\varphi)
\right) = 0 .$$ 
In the special one dimensional case, the previous equation becomes
$$
\frac{d}{d\varphi} 
\left(e^{kH(\I+{\partial u_k\over \partial \varphi},\varphi)} 
\frac{\partial H}{\partial  I}(\I+{\partial u_k\over \partial \varphi},\varphi)
\right) = 0.
$$
This can be integrated and one obtains
\begin{equation}
e^{kH(\I+{\partial u_k\over \partial \varphi} ,\varphi)} 
\left( \frac{\partial H}{\partial  I}(\I+{\partial u_k\over \partial \varphi},\varphi)  \right) = c 
\label{eqc}
\end{equation}
for some $c \in {\Bbb R}$. For $H(I,\varphi)=I^2/2+f(\varphi)$, we have
\begin{equation}
e^{kH(\I+{\partial u_k\over \partial \varphi} ,\varphi)} 
\left( \I+{\partial u_k\over \partial \varphi} \right) = c  .
\label{eqcsimp}
\end{equation}
From equation (\ref{eqcsimp}) one immediately recognizes that the constant $c$ 
has the same sign as $\I+{\partial u_k\over \partial \varphi}$ and 
 $\I$. It suffices to first write (\ref{eqcsimp}) as
$\I+{\partial u_k\over \partial  \varphi}= 
c e^{-k H(\I+{\partial u_k\over \partial \varphi} ,\varphi)}$, and then to average 
both sides over $\varphi$. Therefore, for $\I>0$, $c>0$ one also  has $\I+{\partial u_k\over \partial \varphi}>0$ for any 
$\varphi\in {\Bbb S}^1$, so that we can write equation (\ref{eqcsimp}) in the 
form
$$
 e^{k \Big (\frac{1}{2} (\I + {\partial u_k\over \partial \varphi})^2+
f(\varphi)\Big )+\log (\I + {\partial u_k\over \partial \varphi})}=c.
$$
Thus, on putting  $c_k := \frac{\log c}{k}$, one has
\begin{equation}
\frac{1}{2} (\I + {\partial u_k\over \partial \varphi})^2 + 
\frac{1}{k} \log(\I+ {\partial u_k\over \partial \varphi})+f(\varphi) = c_k 
 .
\label{HJJ}
\end{equation}
Bearing in mind the Lambert function $W$ (see Section \ref{W}), we see that
the equation 
$\frac{1}{2}\gamma_k^2 + \frac{1}{k} \log(\gamma_k) + f(\varphi) = c_k$ may be 
written in the form
$$
e^{k \gamma_k^2}(k \gamma_k^2) = e^{2 k (c_k-f(\varphi ))}k ,
$$
and since the right hand side is positive, we can represent its solution --see formula (\ref{lambertfc})-- as
\begin{equation} \label{GAMMA CAPPA}
k \gamma_k^2= W(e^{2(c_k-f(\varphi))k}k)  ,
\end{equation}
that is,
\begin{equation} 
\I+{\partial u_k\over \partial \varphi} = \gamma_k(c_k,\varphi) = 
\sqrt{\frac{W(e^{2(c_k-f(\varphi))k}k)}{k}} .
\label{ipreintcic}
\end{equation}
Integrating (\ref{ipreintcic}) over $[0,\varphi]$ we find
$$
u_k(\I,\varphi) = u_k(\I,0) - \I \varphi + \int^{\varphi}_0 \gamma_k(c_k,x)dx .
$$
If we now require that $u_k$  be 
periodic with respect to $\varphi$, we have
$$
\I = \frac{1}{2 \pi} \int_0^{2\pi} \gamma_k(c_k,\varphi )d\varphi ,
$$
while a function with zero average is obtained if one chooses
$$
u_k(\I,0) = \I \pi - 
\frac{1}{2\pi} \int^{2\pi}_{0} \int^{\varphi}_{0} \gamma_k(c_k,x) dx d\varphi
 .
$$
We have therefore proved that the function $u_k(\I,\varphi)$ in
(\ref{soluzione}) has zero average and solves the Euler-Lagrange equation
(\ref{divsigma}) for $I_k[u]$. \\
From  definitions  (\ref{sigma}) and (\ref{acca}), and since
$$
e^{k H(\tilde I+{\partial u_k\over \partial \varphi},\varphi)}= 
e^{k(c_k -{1\over k}\log \gamma_k)}={e^{k c_k}\over \gamma_k},
$$
we immediately obtain (\ref{sigmak}) and
(\ref{hamiltoniana}). The limit (\ref{liminthams}) follows 
directly from \cite{E2} (specifically from (2.5) in \cite{E2}, see also 
Theorems 2.1 and 3.1 in \cite{E1}), while structure (\ref{H-lim}) comes from 
the well known representation
of the effective Hamiltonian for one dimensional systems. \hfill $\Box$

\subsection{Uniform convergence of $\gamma_k$ to $\gamma_0$}
This section is devoted to proving the uniform convergence of
$\gamma_k(c,\varphi)$ to $\gamma_0(c,\varphi)$ on compact sets of $\mathbb{R}
\times \mathbb{S}^1$. Specifically, we prove that for any
$\varepsilon>0$ and $c_*\in {\Bbb R}$,  there exists 
$K(\epsilon,c_*)$ such that, for any $k \ge K(\epsilon,c_*)$, we have
\begin{equation*}
\norm{ \gamma_k(c,\varphi)-\gamma_0 (c,\varphi)}\leq \varepsilon
\label{unifgamma}
\end{equation*}
$\forall (c,\varphi) \in (-\infty,c_*] \times {\Bbb S}^1$. This result will be
  essential in the proof of the pointwise convergence of $c_k$ to $c$. We distinguish two different cases $(i)$ and $(ii)$. 
\begin{itemize}
\item[$(i)$] Let us consider $\varphi$ such that $c \geq f(\varphi)$. We start 
with the following
\begin{Lem}\label{lemma1}
Let $c\in {\Bbb R}$. For any $\varepsilon >0$ there exists $K_0(\varepsilon)$ 
independent of $c$ such that, for any $k > K_0(\varepsilon)$ 
and $\varphi$ satisfying $c\geq f(\varphi)$, we have
\begin{equation}
{\norm{\gamma_k(c,\varphi ) - \sqrt{2 (c-f(\varphi)) +{\log k\over k}}}\over 
\sqrt{2 (c-f(\varphi))+{\log k\over k}}}
\leq \varepsilon .
\label{lemmineq}
\end{equation}
\end{Lem}
\textit{Proof.} 
From (\ref{liminf}) we know that for any $\varepsilon > 0$ there exists $K_0(\varepsilon)$ such that, for any $z > K_0(\varepsilon)$, we have
$$\norm{\sqrt{\frac{W(z)}{\log z}} - 1} \le \varepsilon  .$$ 
Moreover, since $c \geq f(\varphi )$ one also has
$$e^{2(c - f(\varphi))k}k \ge k  .$$
As a consequence of the last two facts, for any $\varepsilon > 0$ 
and $k > K_0(\epsilon)$, we have
$$
\norm{ \sqrt{W( e^{2 (c -f(\varphi ))k}k) \over 2 (c-f(\varphi))k+\log k}-1}
\leq \varepsilon
$$
$\forall \varphi$ such that $c\geq f(\varphi)$. We write the above inequality 
as
$$
{\norm{\sqrt{{W( e^{2 (c -f(\varphi ))k}k )\over k}} - \sqrt{2 (c-f(\varphi)) +{\log k\over k}}}\over 
\sqrt{2 (c-f(\varphi))+{\log k\over k}}}
\leq \varepsilon ,
$$
from which the lemma immediately follows. \hfill $\Box$ \\
~\newline
\noindent The uniform convergence  of $\gamma_k(c,\varphi)$ to 
$\gamma_0 (c,\varphi)$ on the set $c\geq f(\varphi)$ 
is now a direct consequence of the lemma. If $c < \min f$  
there is nothing 
to prove, since this set is empty. We can therefore assume $c_*\geq c 
 \geq \min f$. For any $\eta>0$, from Lemma \ref{lemma1}
there exists $K_0(\eta)$ such that, for any $k \geq 
K_0(\eta)$, we have
$$
\norm{ \gamma_k(c,\varphi)-\gamma_0 (c,\varphi)}\leq 
\eta \sqrt{2 (c-f(\varphi))+{\log k\over k}}+
\norm{\sqrt{2 (c-f(\varphi)) +{\log k\over k}}-
\sqrt{2 (c-f(\varphi))}} .
$$ 
Choosing $K_1(\eta)$ such that ${\log k\over k}\leq \eta$ for any 
$k \geq K_1(\eta)$, we immediately obtain 
\begin{equation} \label{uniform estimate}
\norm{ \gamma_k(c,\varphi)-\gamma_0 (c,\varphi)}\leq  
\eta \sqrt{2 (c-\min f) +\eta}+\sqrt{\eta} \le \eta \sqrt{2 (c_* - \min f) +\eta}+\sqrt{\eta}
\end{equation}
$\forall k \geq \max \{ K_0(\eta), K_1(\eta)\}$. Therefore, if for any $\epsilon>0$ we 
choose $\eta:=\eta(\epsilon,c_*)$ such that
$$
\eta(\epsilon,c_*) \sqrt{2 (c_* - \min f) +\eta(\epsilon,c_*)}+
\sqrt{\eta(\epsilon,c_*)} = \epsilon ,
$$
we find that, for any $k\geq  \max \{
K_0(\eta(\epsilon,c_*)), K_1(\eta(\epsilon,c_*))\}$,  one has
$$
\norm{ \gamma_k(c,\varphi)-\gamma_0 (c,\varphi)}\leq  \epsilon.
$$
\item[$(ii)$] We now consider $\varphi$ such that $c < f(\varphi)$.
In this  case $\gamma_0 (c,\varphi)=0$ and therefore
$$
\norm{\gamma_k(c,\varphi)-\gamma_0 (c,\varphi)} = 
\sqrt{W( e^{2 (c -f(\varphi ))k}k) \over k} .
$$
Since $W$ is an increasing function of $z\in [0,+\infty)$ and $e^{2 (c -f(\varphi ))k}k \leq k$, we have
$$
\sqrt{W( e^{2 (c -f(\varphi ))k}k) \over k} \leq 
\sqrt{W(k) \over k} .
$$
From (\ref{liminf}) we obtain
$$
\lim_{k\rightarrow +\infty} \sqrt{ {W(k)\over \log k}}=1, 
$$
that is, for any $\eta>0$ there exists $K_0(\eta)$ such that, for 
every $k \geq K_0(\eta)$, 
$$
\norm{\sqrt{{W(k) \over \log k}} -1} \leq \eta \ \ \Longrightarrow \ \ 
\norm{\sqrt{W(k)} - \sqrt{\log k}} \le \eta \sqrt{\log k}  ,
$$
from which we have
$$
\sqrt{W(k) \over k} \leq {\sqrt{\log k} + \norm{ {\sqrt{W(k)}-\sqrt{\log k}} }
  \over \sqrt{k}} \leq \sqrt{{\log k\over k}}(1 + \eta)  .
$$
Therefore, on choosing  $k\geq K_1(\eta)$, it follows that
\begin{equation}
\norm{ \gamma_k(c,\varphi)-\gamma_0 (c,\varphi)}\leq \sqrt{\eta}(1+\eta) .
\label{cminphi}
\end{equation}
For any $\epsilon>0$, we choose $\eta :=\tilde \eta(\epsilon)$ 
such that $\sqrt{\tilde \eta(\epsilon)}(1+\tilde \eta(\epsilon))=\epsilon$. 
\par\noindent
Thus, 
for any \hbox{$k\geq \max \{K_1(\tilde \eta(\epsilon)), K_0(\tilde \eta(\epsilon))\}$,} we have 
$$
\norm{ \gamma_k(c,\varphi)-\gamma_0 (c,\varphi)}\leq \epsilon .
$$
\end{itemize}
The uniform convergence is therefore proved by choosing
$$
K(\epsilon,c_*):=\max \{K_0(\eta(\epsilon,c_*)), K_1(\eta(\epsilon,c_*)),
K_1(\tilde \eta(\epsilon)),
K_0(\tilde \eta(\epsilon))\} .
$$

\subsection{Pointwise convergence of $c_k$ to $c$}

In this section we prove that, for any $\I>0$, we have
$$
\lim_{k\rightarrow +\infty} c_k(\I)=c(\I) .
$$ 
The proof is structured into points $(i)-(iv)$.
\begin{itemize}
\item[$(i)$] We first establish that the sequence $c_k(\I)$ is bounded from
  above. \\
On the contrary, let us suppose the existence of a diverging sub-sequence $c_{k_i}(\I)$:
$$
\lim_{i\rightarrow \infty} c_{k_i}(\I)=+\infty .
$$
From the monotonicity of $W$, for 
any $\varphi\in {\Bbb  S}^1$, we have
$$
\sqrt{W(e^{2(c_{k_i}(\I)-f(\varphi))k_i}k_i)
\over k_i} \geq  \sqrt{W(e^{2(c_{k_i}(\I)-\max f )k_i}k_i)
\over k_i}
 ,
$$
so that, by integrating in $\varphi\in [0,2\pi]$ and using (\ref{ark}) and (\ref{gkdin}), we obtain
\begin{equation}
\I \geq \sqrt{W(e^{2(c_{k_i}(\I)-\max f )k_i}k_i)
\over k_i} .
\label{confronti}
\end{equation}

Moreover, as a consequence of (\ref{liminf}), the divergence of $c_{k_i}(\I)$ 
implies the divergence of the 
sequence $a_i = \sqrt{W(e^{2(c_{k_i}(\I)-\max f )k_i}k_i)
\over k_i}$. Indeed one has
$$
\lim_{i \to +\infty} a_i^2 = \lim_{i\rightarrow +\infty} 
{W(e^{2(c_{k_i}(\I)- \max f )k_i}k_i)\over 
2(c_{k_i}(\I)-\max f)k_i+\log k_i}\Big (2(c_{k_i}(\I)-\max f)+
{\log k_i\over k_i}\Big )=+\infty .
$$
But this is in contradiction with inequality (\ref{confronti}).

\item[$(ii)$] We proceed by proving that, for $\I>0$, there exists $K_2(\I)$ 
such that 
$$
c_k(\I) > c(\I/4) >\min f 
$$
$\forall k\geq K_2(\I)$.  

Let us first prove that $c_k(\I) > c(\I/4)$ for sufficiently large $k$.
Point $(i)$ provides the existence of $c_*(\I)$ for which 
$\sup_k c_k(\I) < c_*(\I)$. Therefore, from the uniform convergence of 
$\gamma_k$ to $\gamma_0$ in 
$(-\infty,c_*(\I)]\times {\Bbb S}^1$, we find that for any $\varepsilon>0$ 
there exists $K(\varepsilon,c_*(\I))$ such that, for any $k\geq  
K(\varepsilon,c_*(\I))$ and $c< c_*(\I)$, one has
$$
{1\over 2\pi}\norm{\int_0^{2\pi}\gamma_k(c,\varphi)d\varphi -
\int_0^{2\pi}\gamma_0(c,\varphi)d\varphi}\leq \varepsilon .
$$
In particular,
$$
{1\over 2\pi}\norm{\int_0^{2\pi}\gamma_k(c_k(\I),\varphi)d\varphi -
\int_0^{2\pi}\gamma_0(c_k(\I),\varphi)d\varphi}=
\norm{\I- {1\over 2\pi}\int_0^{2\pi}\gamma_0(c_k(\I),\varphi)d\varphi}
\leq \varepsilon .
$$
From the above inequality we immediately obtain
\begin{equation}
\I \leq \varepsilon +{1\over 2\pi}\int_0^{2\pi}
\gamma_0(c_k(\I),\varphi)d\varphi .
\label{Ieps}
\end{equation}

We proceed by considering the function
$$
\I(c)={1\over 2\pi}\int_0^{2\pi}\gamma_0(c,\varphi)d\varphi ,
$$
which is strictly monotone for $c\geq \min f$. Moreover, if we 
fix $\epsilon=\I/2$, the inequality (\ref{Ieps}) gives
\begin{equation}
{\I\over 2} \leq \I(c_k(\I))
\label{ickmezzi}
\end{equation}
$\forall k\geq K(\I/2,c_*(\I))$. Inequality (\ref{ickmezzi}) also implies that
$$
c_k(\I)> c(\I/4)
$$
$\forall k\geq K(\I/2,c_*(\I))$. In fact, if there exists 
$k\geq K(\I/2,c_*(\I))$ such that $c_k(\I)\leq c(\I/4)$, from (\ref{ickmezzi}) 
and the monotonicity of $\I(c)$, we have also
$$
{\I\over 2} \leq \I(c_k(\I)) \leq \I \Big ( c\Big ({\I\over 4}\Big )\Big )
={\I\over 4} ,
$$
which is a contradiction. Moreover, since $\I>0$,  one necessarily has $c(\I/4)>
\min f$.

\item[$(iii)$] We prove that, for any $c''\geq c'> \min f$, we have
\begin{equation}
\norm{\I(c'')-\I(c')} \geq \norm{c''-c'} {m(c')\over 2\sqrt{2}(c''-\min f)} ,
\label{eqicc}
\end{equation}
where $m(c')$ is the measure of the set
$$
A_+(c'):=\{\varphi \in {\Bbb S}^1: c'\geq f(\varphi)\} .
$$
Indeed, since $\I(c)$ is a strictly monotone,
\begin{eqnarray}
\norm{\I(c'')-\I(c')} &=& {1\over 2\pi}\int_0^{2\pi}
(\gamma_0(c'',\varphi)-\gamma_0(c',\varphi))d\varphi \geq\cr 
&\geq& {1\over 2\pi}\int_{A_+(c')} 
(\sqrt{2(c''-f(\varphi))}-
\sqrt{2(c'-f(\varphi))}) d\varphi =\cr
&=&{\norm{c''-c'}\over \pi} \int_{A_+(c')}
{1\over \sqrt{2(c''-f(\varphi))}+\sqrt{2(c'-f(\varphi))}}d\varphi\geq\cr
&\geq& \norm{c''-c'} {m(c')\over 2\sqrt{2}\pi (c''-\min f)} .
\end{eqnarray}

\item[$(iv)$] Finally, from point $(ii)$ we know that for any $\eta>0$, there exists 
$K(\eta,c_*(\I))$ such that,
for any $k \geq K(\eta,c_*(\I))$,
$$
\norm{\I(c(\I))-\I(c_k(\I))}
=\norm{\I- {1\over 2\pi}\int_0^{2\pi}\gamma_0(c_k(\I),\varphi)d\varphi}
\leq \eta .
$$
Moreover, since $c(\I)>\min f$ and  for any 
$k\geq K_2(\I)$ one also has $c_k(\I)>c(I/4)>\min f$, we can apply 
inequality (\ref{eqicc}) to $c(\I)$ and $c_k(\I)$, obtaining
$$
\norm{\I(c(\I))-\I(c_k(\I))} \geq 
\norm{c(\I)-c_k(\I)} {m(\min \{c(\I),c_k(\I)\})\over 
2\sqrt{2}\pi(\max\{c(\I),c_k(\I)\}-\min f)}
$$
and therefore
$$
\norm{c(\I)-c_k(\I)}  \leq 2\sqrt{2}\pi\  \eta\  {(\max\{c(\I),c_k(\I)\}-\min f) 
\over m(\min \{c(\I),c_k(\I)\})} \leq 
 2\sqrt{2}\pi \ \eta\ {\max\{c(\I),c_*\}\over m(c(\I/4))} .
$$
Therefore, setting
$$
\eta = \epsilon {m(c(\I/4)) \over  2\sqrt{2}\pi \max\{c(\I),c_*(\I)\}}
$$
for any $\epsilon>0$, we have proved that there exists $K_3(\epsilon,\I)$ such that
$$
\norm{c(\I)-c_k(\I)} \leq \epsilon .
$$
$\forall k \geq K_3(\epsilon,\I)$.
\end{itemize}
\subsection{Uniform convergence of $u_k$ to $u_0$}

In this section we prove that for any $\I>0$ and $\varepsilon > 0$ there exists $\tilde K(\varepsilon,\I)$ such
that, for any $k \ge \tilde K(\varepsilon,\I)$, we have
$$|u_k(\I,\varphi) - u_0(\I,\varphi)| \le \varepsilon$$
$\forall \varphi \in \mathbb{S}^1$. \\ \\
The result follows from the estimates (I) and (II) to be established below.
\begin{itemize}
\item[(I)] 
We know from $(i)$ of Section 5.3 that there exists $c_{*}(\I)$ such that
$\sup_{k} c_k(\I) < c_{*}(\I)$. We can therefore apply the convergence result
of Section 5.2 to conclude that for any $\eta > 0$ there exists
$K(\eta, c_{*}(\I))$ such that, for any $k \ge K(\eta,c_{*}(\I))$,
one has
\begin{equation} \label{gamma c}
|\gamma_k(c_k(\I),\varphi) - \gamma_0(c_k(\I),\varphi)| \le \eta
\end{equation}
$\forall \varphi \in \mathbb{S}^1$. 
\item[(II)] For any $c\in {\Bbb R}$,
$$
\lim_{c'\rightarrow c}\max_{\varphi \in {\Bbb S}^1} 
\norm{\gamma_0(c',\varphi)-\gamma_0(c,\varphi)}=0 .
$$ 
In other words, for any 
$\eta>0$ there exists $\rho(\eta,c)$ such that for any $c'$ with 
$\norm{c'-c}\leq \rho(\eta)$ and $\varphi\in {\Bbb S}^1$,
\begin{equation}
\norm{\gamma_0(c',\varphi)-\gamma_0(c,\varphi)}\leq \eta .
\label{contgamm}
 \end{equation}
This is trivial if $c\leq \min f$ or $c\geq \max f$. Let us therefore consider
$c\in (\min f,\max f)$ and distinguish two cases. \\
(i) Let 
$\varphi$ be such that $c-f(\varphi)\leq \eta^2/32$, 
so that $\gamma_0(c,\varphi) \leq \eta/4$. If $\norm{c'-c}\leq \eta^2/32$, 
then we have also $c'-f(\varphi)\leq {\eta^2/16}$ and
$$
\norm{\gamma_0(c',\varphi)-\gamma_0(c,\varphi)} \leq 
\gamma_0(c',\varphi)+\gamma_0(c,\varphi) \leq \eta .
$$
(ii) Let  $\varphi$ be such that $c-f(\varphi)\geq \eta^2/32$, so that
$\gamma_0(c,\varphi) =\sqrt{2(c-f(\varphi))}\geq \eta/4$. If $\norm{c'-c}\leq \eta^2/64$,
then we also have $c'-f(\varphi)\geq \eta^2/64$ and
$\gamma_0(c',\varphi) =\sqrt{2(c'-f(\varphi))}\geq \eta/(4\sqrt{2})$. As a consequence
$$
\norm{\gamma_0(c',\varphi)-\gamma_0(c,\varphi)} \leq  
{2\norm{c'-c}\over \sqrt{2(c-f(\varphi))}+\sqrt{2(c'-f(\varphi))}}\leq
{ 2 {\eta^2\over 64}\over {\eta\over 4}+{\eta\over 8}}\leq \eta .
$$
Therefore the uniform continuity is proved also for $c\in (\min f,\max f)$ 
with $\rho(\eta)=\eta^2/64$. 
\end{itemize}
Hence
\begin{eqnarray*}
|u_k(\I,\varphi) - u_0(\I,\varphi)| &=& \norm{ \frac{1}{2 \pi} \int^{2\pi}_0 \int^{\varphi}_0 [\gamma_k(c_k(\I),x) - \gamma_0(c(\I),x)] dx d\varphi + \int^{\varphi}_0 [\gamma_k(c_k(\I),x) - \gamma_0(c(\I),x)] dx} \\
&\le& 4\pi \sup_{x\in {\Bbb S}^1}
\norm{\gamma_k(c_k(\I),x) - \gamma_0(c(\I),x)} \\
&\le& 4\pi \sup_{x\in {\Bbb S}^1}\Big (
|\gamma_k(c_k(\I),x) - \gamma_0(c_k(\I),x)| +
|\gamma_0(c_k(\I),x) - \gamma_0(c(\I),x)|\Big ).
\end{eqnarray*}
From (I), for any $k\geq K(\epsilon/8\pi,c_*(\I))$, we have
$$
4\pi \sup_{x\in {\Bbb S}^1}
|\gamma_k(c_k(\I),x) - \gamma_0(c_k(\I),x)| \leq {\epsilon\over 2} .
$$
Since $c_k(\I)$ converges to $c(\I)$, for any $k \geq K_3(\rho(\epsilon/8\pi),
\I)$ we have $\norm{c_k(\I)-
c(\I)}\leq \rho(\epsilon/8\pi)$, and therefore by (II)
$$
  4\pi \sup_{x\in {\Bbb S}^1} |\gamma_0(c_k(\I),x) - \gamma_0(c(\I),x)|
\leq {\epsilon \over 2} .
$$
For $k \geq \tilde K(\epsilon,\I)= 
\max \{  K(\epsilon/8\pi,c_*(\I)), K_3(\rho(\epsilon/8\pi),\I )\}$, we can 
therefore write
\begin{eqnarray*}
|u_k(\I,\varphi) - u_0(\I,\varphi)| \le  \varepsilon 
\end{eqnarray*}
$\forall \varphi \in \mathbb{S}^1$, which establishes the statement on uniform convergence.
\hfill $\Box$

\section{Proof of Theorem \ref{prop2}} \label{sei}
Before treating the one dimensional case, we show that for general mechanical Hamiltonian systems with  $n$ degrees of freedom,
\begin{equation} \label{MECH}
H(I,\varphi) = \sum_{i=1}^n\frac{I_i^2}{2} + f(\varphi_1,\ldots,\varphi_n),
\end{equation}
the integrals in (\ref{estsigmakeff}) and (\ref{estsigmak}) can be respectively written in the form
\begin{equation} \label{int1}
E_1(k) =  
\int_{|\tilde{I}| \le R} \int_{\mathbb{T}^n} \Big \vert
\frac{\partial}{\partial I} \Big[\frac{1}{2}(\tilde{I} + \frac{\partial
    u_k}{\partial \varphi})^2 - \bar{H}_k(\tilde{I})\Big] \Big\vert^2 \sigma_k
 d\varphi d\tilde{I}
\end{equation}
and
\begin{equation} \label{int2}
E_2(k) = 
\int_{\mathbb{T}^n} \Big \vert \frac{\partial}{\partial \varphi} \Big[
\frac{1}{2}(\tilde{I} + \frac{\partial u_k}{\partial \varphi})^2 +
f(\varphi)\Big] \Big\vert^2  \sigma_kd\varphi  .
\end{equation}
\indent Indeed, from (\ref{genfun2}), and with the notation of Section \ref{SEZIONE 2}, we have
$$\dot{\tilde{\varphi}}^t = \dot{\varphi}^t + \frac{\partial^2 u_k}{\partial \varphi \partial I}(\tilde{I},\varphi^t) \dot{\varphi}^t  .$$
Therefore, since for the Hamiltonian (\ref{MECH}) the torus flow is
$$\dot{\varphi}^t = \tilde{I} + \frac{\partial u_k}{\partial \varphi}(\tilde{I},\varphi^t)  ,$$
we have
$$
\dot{\tilde{\varphi}}^t = (\mathbb{I} + \frac{\partial^2 u_k}{\partial \varphi \partial I})(\tilde{I} + \frac{\partial u_k}{\partial \varphi}) = \frac{\partial}{\partial I} \Big[\frac{1}{2}(\tilde{I} + \frac{\partial u_k}{\partial \varphi})^2\Big]_{\vert \varphi=
{\cal C}^t_{\tilde I}(\varphi)} 
$$
from which we immediately obtain
$$
\int_{|{\tilde I}| \leq R} \int_{{\Bbb T}^n} \Big \vert {\partial {\tilde \varphi}^t \over \partial t} -
{\partial {\bar H}_k \over \partial I} (\tilde I) \Big \vert^2 \sigma_k d\varphi d\tilde I = 
\int_{|\tilde{I}| \le R} \int_{\mathbb{T}^n} \Big \vert
\frac{\partial}{\partial I} \Big[\frac{1}{2}(\tilde{I} + \frac{\partial
    u_k}{\partial \varphi})^2 - \bar{H}_k(\tilde{I})\Big]_{\vert \varphi=
{\cal C}^t_{\tilde I}(\varphi)} \Big\vert^2 \sigma_k
 d\varphi d\tilde{I} .
$$
Since the torus flow preserves the measure defined by $\sigma_k$, 
see (\ref{divsigma}), we obtain (\ref{int1}). \\
\indent To prove (\ref{int2}), on recalling  (\ref{genfun1}):
$$I^t = \tilde{I} + \frac{\partial u_k}{\partial \varphi}(\tilde{I},\varphi^t),$$
we now obtain
$$\dot{I}^t = \frac{\partial^2 u_k}{\partial \varphi^2}(\tilde{I},\varphi^t) \dot{\varphi}^t = \frac{\partial^2 u_k}{\partial \varphi^2}(\tilde{I},\varphi^t) (\tilde{I} + \frac{\partial u_k}{\partial \varphi}(\tilde{I},\varphi^t))  .$$
Moreover, for mechanical Hamiltonians (\ref{MECH}),
$$\frac{\partial H}{\partial \varphi}(\Phi^t(\tilde{I},\varphi)) = \frac{\partial f}{\partial \varphi}(\varphi^t)  ,$$
hence
$$\dot{I}^t + \frac{\partial H}{\partial \varphi}(\Phi^t(\tilde{I},\varphi)) =
\frac{\partial^2 u_k}{\partial \varphi^2}(\tilde{I},\varphi^t) (\tilde{I} +
\frac{\partial u_k}{\partial \varphi}(\tilde{I},\varphi^t)) + {\partial f
\over \partial \varphi}(\varphi^t) = 
\frac{\partial}{\partial \varphi} \Big[ 
\frac{1}{2}(\tilde{I} + \frac{\partial u_k}{\partial \varphi})^2 + 
f(\varphi)\Big]_{\vert \varphi=
{\cal C}^t_{\tilde I}(\varphi)}
$$
and, again using (\ref{divsigma}), formula (\ref{int2}) is proved. 

\subsection{The one dimensional case}
In the one dimensional case, for $\I>0$ one has  
$\gamma_{k} = \tilde{I} + \frac{\partial u_k}{\partial \varphi}$.
Thus on invoking  the symmetry with respect to $\tilde I$, 
formulas (\ref{int1}) and (\ref{int2}) become respectively
\begin{equation} \label{inte1}
E_1(k)=\int^R_{-R} \int^{2\pi}_0 \Big \vert {\partial {\tilde \varphi}^t \over \partial t} -
{\partial {\bar H}_k \over \partial I} (\tilde I) \Big \vert^2 \sigma_k
d\varphi d\tilde I = 2 
\int^R_0 \int^{2\pi}_0 \Big \vert
\frac{\partial}{\partial I} \Big[\frac{\gamma_k^2}{2} - \bar{H}_k(\tilde{I})\Big] \Big\vert^2 \sigma_k d\varphi d\tilde{I}
\end{equation}
and 
\begin{equation} \label{inte2}
E_2(k)=\int^{2\pi}_0 \Big\vert 
\frac{\partial I^t}{\partial t} + 
\frac{\partial H}{\partial \varphi}(\Phi^t(\tilde{I},\varphi))
\Big \vert^2d\varphi = \int^{2\pi}_0 \Big \vert \frac{\partial}{\partial \varphi} \Big[
\frac{\gamma_k^2}{2} + f(\varphi)\Big] \Big\vert^2\sigma_kd\varphi  .
\end{equation}

\subsection{Proof of (\ref{stima1}) for $E_1(k)$}

From (\ref{HJ}) and (\ref{hamiltoniana}) we obtain
\begin{eqnarray} 
E_1(k)&=& 2 
\int^R_0 \int^{2\pi}_0 \Big \vert
\frac{\partial}{\partial I} \Big[c_k(\I)-{1\over k}
\log \gamma_k(c_k(\I),\varphi) - 
\bar{H}_k(\tilde{I})\Big] \Big\vert^2 \sigma_k d\varphi d\tilde{I}
=\cr
&=& 2\int^R_0 \int^{2\pi}_0 \Big \vert
\frac{\partial}{\partial I} \Big[-{1\over k}
\log \gamma_k(c_k(\I),\varphi) - 
{1\over k}\log \int_0^{2\pi}{dx\over \gamma_k(c_k(\I),x)}
] \Big\vert^2 \sigma_k d\varphi d\tilde{I} ,
\end{eqnarray}
that is
$$
E_1(k)={2\over k^2}
\int^R_0 \int^{2\pi}_0 \Big \vert
{ \gamma_k'(c_k,\varphi)\over \gamma_k(c_k,\varphi)}-
{\int_0^{2\pi}{\gamma_k'(c_k,x)\over \gamma_k^2(c_k,x)}dx
\over \int_0^{2\pi}{dx\over \gamma_k(c_k,x)}}
\Big\vert^2 \sigma_k d\varphi d\tilde{I} ,
$$
where we have denoted
$$
\gamma_k'(c_k,\varphi):={\partial \over \partial \I}
\gamma_k(c_k(\I),\varphi),\qquad \gamma_k'(c_k,x):=
{\partial \over \partial \I}\gamma_k(c_k(\I),x) .
$$
Now, by expanding the square, using (\ref{sigmak}) and 
performing the integration over $\varphi$, we have
$$
E_1(k)={2\over k^2}
\int^R_0 {\int_0^{2\pi}{{\gamma'_k}^{2}(c_k,x)\over \gamma_k^3(c_k,x)}dx
\int_0^{2\pi}{dx\over \gamma_k(c_k,x)} - \Big (
\int_0^{2\pi}{\gamma_k'(c_k,x)\over \gamma_k^2(c_k,x)}dx \Big )^2 
\over \Big ( \int_0^{2\pi}{dx\over \gamma_k(c_k,x)} \Big )^2}d\tilde{I} .
$$
From (\ref{HJ}) it is also easy to see that
\begin{equation*} \label{gamma-c}
\frac{d \gamma_k}{d c}\Big{|}_{c=c_k} = \frac{k\gamma_k}{k\gamma_k^2 + 1} ,
\end{equation*} 
whence
\begin{equation} \label{gamma-I}
\gamma'_k(c_k,x)=
\frac{d \gamma_k}{d c}\Big{|}_{c=c_k} c'_k
= \frac{\gamma_k(c_k,x)}{{1\over k}+\gamma_k^2(c_k,x)} c'_k .
\end{equation} 
Using (\ref{gamma-I}) we have the representation
\begin{equation}
E_1(k)={2\over k^2} \int^R_0 {c'_k}^2(\I)
{A_{-{1\over 2},2}(k,c_k(\I))A_{-{1\over 2},0}(k,c_k(\I))
-A_{-{1\over 2},1}^2(k,c_k(\I))\over A_{-{1\over 2},0}^2(k,c_k(\I))}
d\I,
\end{equation}
where the functions
\begin{equation*}
A_{\alpha,\beta}(k,c) := \int_0^{2\pi}\frac{\gamma_k^{2\alpha}(c,\varphi)}{\Big ( {1\over k}+\gamma_k^2(c,\varphi)\Big
    )^\beta}d\varphi
\label{Aalpha}
\end{equation*}
are defined for any $c\in {\Bbb R}$. \\
Since $R>\I(\max f+r)>0$, we will proceed with the integral
\begin{equation}
\tilde E_1(k)={2\over k^2} \int^R_{\I(\max f+r)} {c'_k}^2(\I)
{A_{-{1\over 2},2}(k,c_k(\I))A_{-{1\over 2},0}(k,c_k(\I))
-A_{-{1\over 2},1}^2(k,c_k(\I))\over A_{-{1\over 2},0}^2(k,c_k(\I))}
d\I ,
\label{e1tild}
\end{equation}
since obviously $E_1(k) \geq \tilde E_1(k)$. From (\ref{actik}) it is 
also easy to obtain
\begin{equation} \label{c-P}
c'_k(\I) = {2\pi \over \int^{2\pi}_0 \frac{k \gamma_k}{k\gamma_k^2 + 1}}
={2\pi \over A_{{1\over 2},1}(k,c_k(\I))} . 
\end{equation}
Since $c'_k(\I)>0$ for any $\I>0$, in particular for $R>\I(\max f+r)>0$,  
we can change the integration
variable in $\tilde E_1(k)$ from $\tilde I$ to $c = c_k(\tilde I)$ 
and then invoking (\ref{c-P}) we have
\begin{equation}
\tilde E_1(k)={4\pi\over k^2}  \int^{c_k(R)}_{c_k(\I(\max f+r))}
{A_{-{1\over 2},2}(k,c)A_{-{1\over 2},0}(k,c)
-A_{-{1\over 2},1}^2(k,c)\over A_{-{1\over 2},0}^2(k,c)A_{{1\over
    2},1}(k,c) }dc .
\label{e1kob}
\end{equation}
In order to estimate the last integral by the Dominated 
Convergence Theorem, we organize the proof into points (i) - (viii).

\begin{itemize}
\item[(i)] Since
$$
\lim_{k\rightarrow +\infty} c_k(R)=c(R), \qquad 
\lim_{k\rightarrow +\infty} c_k(\I(\max f+r))=c(\I(\max f+r))=\max f+r
$$
there exist $K_4(r),K_5(R)$ such that, for any $k\geq K_4(r)$,
\begin{equation}
c_k(\I(\max f+r)) > \max f+{r\over 2}
\label{clow}
\end{equation}
and for any $k\geq K_5(R)$,
\begin{equation}
c_k(R)< 2 c(R) .
\label{cip}
\end{equation}
\item[(ii)] For any $\epsilon >0$, $R>\I(\max f+r)$,  
there exists $K_6(\epsilon,R)$ such that, for any $k\geq
  K_6(\epsilon,R)$ and $c\in [c_k(\I(\max f+r)),c_k(R)]$,
\begin{equation}
\norm{\gamma_k(c,\varphi)-\gamma_0(c,\varphi)}<\epsilon
\end{equation}
$\forall \varphi$. This follows from (\ref{cip}) and 
the uniform convergence of $\gamma_k$ to $\gamma_0$ on the set $c\leq  2
c(R)$, as soon as $K_6(\epsilon,R)=\max \{ K_5(R),K(\epsilon,2 c(R))\}$.  
\item[(iii)] For any $k\geq \max \{ K_4(r),K_6(\sqrt{r}/2,R) \}$, for any 
$c\in [c_k(\I(\max f+r)),c_k(R)]$ and any $\beta \geq 0$, we have
\begin{equation}
A_{-{1\over 2},\beta}(k,c)\leq 2\pi \Big ({4\over r}\Big )^{1+2\beta \over 2}
 .
\end{equation}
Indeed, since $k\geq K_6(\sqrt{r}/2,R)$, one has
$$
\norm{\gamma_k(c,\varphi)-\gamma_0(c,\varphi)}< {\sqrt{r}\over 2}
$$
so that
$$
\gamma_k(c,\varphi) \geq \gamma_0(c,\varphi)-{\sqrt{r}\over 2}
$$
$\forall \varphi$. Since $c\geq c_k(\I(\max f+r))$, if $k\geq K_4(r)$, we have also $c\geq c_k(\I(\max f+r)) \geq \max f+r/2$, and therefore
$$
\gamma_0(c,\varphi)= \sqrt{2(c-\max f)}\geq \sqrt{r} ,
$$
hence
$$
\gamma_k(c,\varphi) \geq {\sqrt{r}\over 2} .
$$
As a consequence,
$$
A_{-{1\over 2},\beta}(k,c)\leq \int_0^{2\pi} {1\over 
\gamma_k(c,\varphi)^{1+2\beta}}d\varphi \leq 2\pi \Big ({4\over r}\Big
)^{1+2\beta\over 2} .
$$
\item[(iv)] For any $k\geq \max \{ K_4(r), K_5(R), K_6(\sqrt{r}/2,R) \}$ and $c\in [c_k(\I(\max f+r)),c_k(R)]$, we have
\begin{equation}
A_{-{1\over 2},0}(k,c)\geq  {2\pi \over  2\sqrt{c(R)}+{\sqrt{r}\over 2}}
 .
\end{equation}
In fact, since $k\geq  K_6(\sqrt{r}/2,R)$, 
$$
\gamma_k(c,\varphi)\leq \gamma_0(c,\varphi) +{\sqrt{r}\over 2} ,
$$
and since $k\geq K_4(r)$ and $k\geq K_5(R)$,
$$
 \gamma_0(c,\varphi)=\sqrt{2(c-f(\varphi))}\leq 
\sqrt{2(2c(R))} .
$$
Therefore
$$
\gamma_k(c,\varphi)\leq  2\sqrt{c(R)}+{\sqrt{r}\over 2} .
$$

\item[(v)]  For any $k\geq \max \{ K_4(r), K_5(R), K_6(\sqrt{r}/2,R) \}$ and $c\in [c_k(\I(\max f+r)),c_k(R)]$, we have
\begin{equation}
A_{{1\over 2},1}(k,c)\geq  \pi{\sqrt{r}\over 1+
(2\sqrt{c(R)}+{\sqrt{r}\over 2})^2} .
\end{equation}
The proof is similar to points (iii) and (iv). 
\item[(vi)] From (iii), (iv) and (v) it immediately follows that 
for any $k\geq \max \{ K_4(r), K_5(R), K_6(\sqrt{r}/2,R) \}$ the integrand in 
(\ref{e1kob}):
$$
{A_{-{1\over 2},2}(k,c)A_{-{1\over 2},0}(k,c)
-A_{-{1\over 2},1}^2(k,c)\over A_{-{1\over 2},0}^2(k,c)A_{{1\over
    2},1}(k,c) }
$$
is dominated by a constant independent of $k,c$.

\item[(vii)] To compute the pointwise limit of the integrand we consider 
the following:
\begin{Lem} \label{lemma62}
Let $c > \max f$, $\lim_{k \to +\infty} c_k = c$ and $\beta\geq \alpha$. Then
\begin{equation} \label{ultimissime}
\lim_{k \to + \infty} A_{\alpha,\beta}(k,c) = a_{\beta - \alpha}(c)  ,
\end{equation}
where
\begin{equation*}
{a}_{\delta}(c) := \int^{2\pi}_0 \frac{1}{\gamma_0^{2\delta}(c,\varphi)} d\varphi = 
\int_0^{2\pi} {1\over [2(c - f(\varphi))]^{\delta}} d\varphi .
\end{equation*}
\end{Lem}
\textit{Proof.} This is an application of the 
Dominated Convergence Theorem. Indeed
$$\lim_{k \to + \infty} \frac{\gamma_k^{2\alpha}(c,\varphi)}{\left( \frac{1}{k} + \gamma_k^2(c,\varphi)\right)^{\beta}} = \frac{1}{\gamma_0^{2(\beta-\alpha)}(c,\varphi)}  .$$
Moreover, the integrand in $A_{\alpha,\beta}(k,c)$:
$$
\frac{\gamma_k^{2\alpha}(c,\varphi)}{\left( \frac{1}{k} +
    \gamma_k^2(c,\varphi)\right)^{\beta}}\leq 
{1\over \gamma_k^{2(\beta-\alpha)}(c,\varphi)}
$$
is dominated by a constant independent on $k,\varphi$ for any $k \geq 
\max \{ K_4(r), K_6(\sqrt{r}/2,R)\}$, see (iii). The statement now follows  
immediately. 
\item[(viii)] Finally, in order to apply the Dominated Convergence Theorem to the 
integral in (\ref{e1kob}), we first write it in the following form
\begin{eqnarray*}
&  \int^{c_k(R)}_{c_k(\I(\max f+r))}
{A_{-{1\over 2},2}(k,c)A_{-{1\over 2},0}(k,c)
-A_{-{1\over 2},1}^2(k,c)\over A_{-{1\over 2},0}^2(k,c)A_{{1\over
    2},1}(k,c) }dc =&  \cr
&=  \int_{\mathbb{R}}  {A_{-{1\over 2},2}(k,c)A_{-{1\over 2},0}(k,c)
-A_{-{1\over 2},1}^2(k,c)\over A_{-{1\over 2},0}^2(k,c)A_{{1\over
    2},1}(k,c) }\chi_{\left[c_k(\I(\max f+r)),c_k(R)\right]}(c) dc &
\end{eqnarray*}
and then we compute its pointwise limit (see Lemma \ref{lemma62}).
$$\lim_{k \to +\infty}   {A_{-{1\over 2},2}(k,c)A_{-{1\over 2},0}(k,c)
-A_{-{1\over 2},1}^2(k,c)\over A_{-{1\over 2},0}^2(k,c)A_{{1\over
    2},1}(k,c) }\chi_{\left[c_k(\I(\max f+r)),c_k(R)\right]}(c)$$
$$= \frac{-a^2_{\frac{3}{2}}(c) + a_{\frac{5}{2}}(c)
 a_{\frac{1}{2}}(c)}{a^3_{\frac{1}{2}}(c)} \chi_{\left[\max f+r,c(R)
 \right]}(c)   .$$
\end{itemize}
Consequently, defining the constant $c_R$ as
\begin{equation*} 
c_R := 4 \pi \int_{\max f + r}^{c(R)} \frac{- a^2_{\frac{3}{2}}(c) + a_{\frac{5}{2}}(c) a_{\frac{1}{2}}(c)}{a^3_{\frac{1}{2}}(c)} dc  ,
\end{equation*}
the estimate (\ref{stima1}) follows. We observe that $c_R \ne 0$ for any non
constant function $f$. Indeed, by using 
the $L^2$-H$\ddot{\text{o}}$lder inequality we obtain
\begin{eqnarray*}
a_{\frac{5}{2}}(c) a_{\frac{1}{2}(c)} &=& \int^{2\pi}_0 \left[\frac{1}{[2(c-f(\varphi))]^{5/4}}\right]^2 d\varphi \int^{2\pi}_0 \left[\frac{1}{[2(c-f(\varphi))]^{1/4}}\right]^2 d\varphi \\
&\ge& \left[\int^{2\pi}_0 \frac{1}{[2(c-f(\varphi))]^{3/2}}d\varphi \right]^2 = a_{3 \over 2}^2(c)  ,
\end{eqnarray*}
hence $c_R \ge 0$, where the equality holds only if $f'(\varphi)=0$ for any 
$\varphi$.  \hfill $\Box$

\subsection{Proof of (\ref{stima2}) for $E_2(k)$}

We know from (\ref{inte2}) that 
\begin{eqnarray*}
E_2(k) &=& \int^{2\pi}_0 \Big{\vert} \frac{\partial}{\partial \varphi} \left[ \frac{\gamma_k^2}{2}(c_k(\tilde{I}),\varphi) + f(\varphi) \right] \Big{\vert}^2\sigma_k(\tilde{I},\varphi) d\varphi \\
&=& \int^{2\pi}_0 \Big{\vert} \left[ \gamma_k(c_k(\tilde{I}),\varphi) \frac{\partial}{\partial \varphi}  \gamma_k(c_k(\tilde{I}),\varphi)
+ f'(\varphi) \right] \Big{\vert}^2\sigma_k(\tilde{I},\varphi) d\varphi .
\end{eqnarray*}
By using the explicit expression for $\gamma_k(c_k(\tilde{I}),\varphi) = \sqrt{\frac{W(e^{2(c_k(\tilde{I}) - f(\varphi))k}k}{k}}$ and the derivative of the Lambert function
$$W'(z) = \frac{1}{(1+W(z))e^{W(z)}}  ,$$
we first establish that
$$E_2(k) = \int^{2\pi}_0 \Big{\vert} \frac{f'(\varphi)}{1 + k \gamma_k^2(c_k(\tilde{I}),\varphi)} \Big{\vert}^2 \sigma_k(c_k(\tilde{I}),\varphi) d\varphi .$$
Moreover, since $$\sigma_k(c_k(\tilde{I}),\varphi) = \frac{1}{A_{-\frac{1}{2},0}(k,c_k(\tilde{I}))} \frac{1}{\gamma_k(c_k(\tilde{I}),\varphi)}  ,$$
we finally obtain
\begin{equation} \label{EDUE}
E_2(k) = \frac{1}{k^2} \frac{1}{A_{-\frac{1}{2},0}(k,c_k(\tilde{I}))} \int^{2\pi}_0 \frac{|f'(\varphi)|^2}{\gamma_k(c_k(\tilde{I}),\varphi)\left(\frac{1}{k} + \gamma_k^2(c_k(\tilde{I}),\varphi)\right)^2}d\varphi  .
\end{equation}
In order to estimate $\lim_{k \to +\infty} k^2 E_2(k)$ by the Dominated Convergence Theorem, we proceed as follows. 
\begin{itemize}
\item[(i)] Since $\sup_{k} c_k(\tilde{I}) < c_{*}(\tilde{I})$, for any $\varepsilon > 0$ there exists $K(\frac{\varepsilon}{2},c_{*}(\tilde{I}))$ such that
$$|\gamma_k(c_k(\tilde{I}),\varphi) - \gamma_0(c_k(\tilde{I}),\varphi)| \le \frac{\varepsilon}{2}$$
$\forall k > K(\frac{\varepsilon}{2},c_{*}(\tilde{I}))$ and $\varphi \in \mathbb{S}^1$. In particular
$$\gamma_k(c_k(\tilde{I}),\varphi) \ge \gamma_0(c_k(\tilde{I}),\varphi) - \frac{\varepsilon}{2}$$
$\forall k > K(\frac{\varepsilon}{2},c_{*}(\tilde{I}))$ and $\varphi \in \mathbb{S}^1$. \\
\item[(ii)]  Moreover, we know that for any $\eta > 0$ there exists $\rho(\eta) > 0$ such that, for any $c'$ with $|c' - c| \le \rho(\eta)$, we have
$$|\gamma_0(c',\varphi) - \gamma_0(c,\varphi)| \le \eta$$
$\forall \varphi \in \mathbb{S}^1$. Therefore, since $\lim_{k \to +\infty} c_k(\tilde{I}) = c(\tilde{I})$, we have
$$|c_k(\tilde{I}) - c(\tilde{I})| \le \rho(\frac{\varepsilon}{2})$$
$\forall k > K_3(\rho(\frac{\varepsilon}{2}),\tilde{I})$. Hence
$$|\gamma_0(c_k(\tilde{I}),\varphi) - \gamma_0(c(\tilde{I}),\varphi)| \le \frac{\varepsilon}{2}$$
$\forall k > K_3(\rho(\frac{\varepsilon}{2}),\tilde{I})$ and $\varphi \in \mathbb{S}^1$. As a consequence,
$$\gamma_0(c_k(\tilde{I}),\varphi) \ge \gamma_0(c(\tilde{I}),\varphi) - \frac{\varepsilon}{2}$$
$\forall k > K_3(\rho(\frac{\varepsilon}{2}),\tilde{I})$ and $\varphi \in \mathbb{S}^1$. \\
\item[(iii)]  Finally, since $c(\tilde{I}) > \max f$ and $\lim_{k \to +\infty} c_k(\tilde{I}) = c(\tilde{I})$, for any $\tilde{I} > 0$ there exists $k(\tilde{I})$ such that
$$c_k(\tilde{I}) > \max f$$ 
$\forall k > k(\tilde{I}).$   
\end{itemize} 
\noindent By using (i) and (ii), we conclude that for any $k > \max \{ K(\frac{\varepsilon}{2},c_{*}(\tilde{I})),K_3(\rho(\frac{\varepsilon}{2}),\tilde{I})\}$ and $\varphi \in \mathbb{S}^1$, we have
$$|\gamma_k(c_k(\tilde{I}),\varphi) - \gamma_0(c(\tilde{I}),\varphi| \le |\gamma_k(c_k(\tilde{I}),\varphi) - \gamma_0(c_k(\tilde{I}),\varphi| + |\gamma_0(c_k(\tilde{I}),\varphi) - \gamma_0(c(\tilde{I}),\varphi| \le \varepsilon.$$
In particular,
\begin{equation} \label{E DUE 1}
\lim_{k \to +\infty} \gamma_k(c_k(\tilde{I}),\varphi) = \gamma_0(c(\tilde{I}),\varphi)
\end{equation}
$\forall \varphi \in \mathbb{S}^1$. \noindent Moreover, for any $k > \max \{ K(\frac{\varepsilon}{2},c_{*}(\tilde{I})),K_3(\rho(\frac{\varepsilon}{2}),\tilde{I}), k(\tilde{I})\}$ and $\varphi \in \mathbb{S}^1$,
\begin{equation} \label{E DUE 2}
\gamma_k(c_k(\tilde{I}),\varphi) \ge \gamma_0(c_k(\tilde{I}),\varphi) - \frac{\varepsilon}{2} \ge \gamma_0(c(\tilde{I}),\varphi) - \varepsilon \ge \sqrt{2(c(\tilde{I}) - \max f)} - \varepsilon  .
\end{equation}
In the last inequalities, since $k > k(\tilde{I})$, the function $\gamma_0(c_k(\tilde{I}),\varphi) = \sqrt{2(c_k(\tilde{I}) - f(\varphi))} > 0$ for any $\varphi \in \mathbb{S}^1$ and therefore
the constant $\sqrt{2(c(\tilde{I}) - \max f)} - \varepsilon$ (independent of $\varphi$ and $k$) is positive (choose, for example, $\varepsilon = \frac{\sqrt{2(c(\tilde{I}) - \max f)}}{2}$). \\
Finally we apply the Dominated Convergence Theorem to compute $\lim_{k \to +\infty} k^2 E_2(k)$, where $E_2(k)$ is given by the expression 
in (\ref{EDUE}). \\
From (\ref{E DUE 2}), we immediately obtain
$$\frac{|f'(\varphi)|^2}{\gamma_k(c_k(\tilde{I},\varphi))\left(\frac{1}{k} + \gamma_k^2(c_k(\tilde{I},\varphi))\right)^2} \le \frac{\max_{\varphi \in \mathbb{S}^1}|f'(\varphi)|^2}{\gamma_k^5(c_k(\tilde{I}),\varphi))} \le \frac{\max_{\varphi \in \mathbb{S}^1}|f'(\varphi)|^2}{\left(\sqrt{2(c(\tilde{I}) - \max f)} - \varepsilon\right)^5} .$$
Moreover, from (\ref{E DUE 1}),
$$\lim_{k \to +\infty} \frac{|f'(\varphi)|^2}{\gamma_k(c_k(\tilde{I}),\varphi)\left(\frac{1}{k} + \gamma_k^2(c_k(\tilde{I}),\varphi)\right)^2} = \frac{|f'(\varphi)|^2}{\gamma_0^5(c(\tilde{I}),\varphi)}  .$$
Therefore
\begin{equation} \label{primo pezzo} 
\lim_{k \to +\infty} \int^{2\pi}_0 \frac{|f'(\varphi)|^2}{\gamma_k(c_k(\tilde{I}),\varphi)\left(\frac{1}{k} + \gamma_k^2(c_k(\tilde{I}),\varphi)\right)^2}d\varphi = \int^{2\pi}_0 \frac{|f'(\varphi)|^2}{\gamma_0^5(c(\tilde{I}), \varphi)}d\varphi .
\end{equation}
We will conclude by proving that
\begin{equation} \label{secondo pezzo}
\lim_{k \to +\infty} A_{-\frac{1}{2},0}(k,c_k(\tilde{I})) = a_{\frac{1}{2}}(c(\tilde{I}))  ,
\end{equation}
where 
$$A_{-\frac{1}{2},0}(k,c_k(\tilde{I})) = \int^{2\pi}_{0} \frac{1}{\gamma_k(c_k(\tilde{I}),\varphi)}d\varphi \quad \text{and} \quad a_{\frac{1}{2}}(c(\tilde{I})) = \int^{2\pi}_0 \frac{1}{[2(c(\tilde{I}) - f(\varphi))]^{1/2}} d \varphi  .$$ 
This limit is  again a straightforward consequence of the Dominated Convergence Theorem. In fact (see (\ref{E DUE 2}) and (\ref{E DUE 1})),
$$
\frac{1}{\gamma_k(c_k(\tilde{I}),\varphi)} \le \frac{1}{\sqrt{2(c(\tilde{I}) - \max f)} - \varepsilon},
$$
where the right hand member is independent of $\varphi$ and $k$, and 
$$\lim_{k \to +\infty} \frac{1}{\gamma_k(c_k(\tilde{I}),\varphi)} = \frac{1}{\gamma_0(c(\tilde{I}),\varphi)}  .$$
From (\ref{primo pezzo}) and (\ref{secondo pezzo}) we therefore obtain the result
$$\lim_{k \to +\infty} k^2 E_2(k) = \frac{1}{a_{\frac{1}{2}}(c(\tilde{I}))} \int^{2\pi}_0 \frac{|f'(\varphi)|^2}{\gamma_0^5(c(\tilde{I}),\varphi)} d\varphi = \frac{1}{a_{\frac{1}{2}}(c(\tilde{I}))} \int^{2\pi}_0 \frac{|f'(\varphi)|^2}{[2(c(\tilde{I}) - f(\varphi)]^{5/2}}d\varphi  .$$

\section{Proof of Proposition \ref{sigma-lim}} \label{sette}
We start by introducing some notation:  $x_k(t)$ and $x(t)$ denote respectively the solutions of
\begin{equation}
\dot x_k=\gamma_k(c_k(\I),x_k) \qquad\textrm{and}\quad \dot x=\gamma_0(c(\I),x),
\label{ficxk}
\end{equation}
 with $x_k(0)=x(0)=0$. 
The function $x_k(t)$ is periodic, with period
$$
T_k := \int_0^{2\pi}{1 \over \gamma_k(c_k(\I),x)} dx  .
$$
\begin{Lem} \label{lemlim1}
If $c(\I)=\max f$, we have 
$$
\lim_{k\rightarrow \infty} T_k=+\infty .
$$
\end{Lem}

\noindent 
\textit{Proof.} We first consider the trivial estimate
$$
\int_0^{2\pi}{1 \over \gamma_k(c_k(\I),x)} dx \geq 
\int_{f(x)\leq c_k(\I)}{1 \over \gamma_k(c_k(\I),x)}dx ,
$$
and then we use Lemma \ref{lemma1} to estimate the right hand side 
integral. In fact, for any $x$ such that $f(x) \leq c_k(\I)$, it follows from 
Lemma \ref{lemma1} that there exist $0<a_1<1<a_2$ and $\tilde K_1(a_1,a_2)$ 
such that
$$
{a_1\over \tilde \gamma_0(c_k(\I),x)} \leq 
{1 \over \gamma_k(c_k(\I),x)} 
\leq {a_2 \over \tilde \gamma_0(c_k(\I),x)}  
$$
for any $k> \tilde K_1(a_1,a_2)$, with 
$\tilde \gamma_0(c_k(\I),x) := \sqrt{2(c_k(\I) - f(x)) + {\log k \over
    k}}$. Therefore, for such $k$ we also have
$$
T_k \geq a_1 \int_{f(x)\leq c_k(\I)} {1 \over \sqrt{2(c_k(\I) - f(x))+
{\log k \over k}}} dx  .
$$
Moreover, from Theorem \ref{prop1}, point $(iii)$, 
for any $\varepsilon > 0$ there exists $\tilde
K_2(\I,\epsilon)$ such that for any $k>\tilde K_2(\I,\epsilon)$} 
it holds $c_k(\I) \leq c(\I) + \varepsilon$. Therefore,
for any $k > \max \{ \tilde K_1(a_1,a_2),\tilde K_2(\I,\epsilon)\}$
we have also
$$
T_k \ge a_1 \int_0^{2\pi} 
{1 \over \sqrt{2(c(\I) + \varepsilon -f(x))+{\log k \over k}}}\chi_{f(x)\leq c_k(\I)}(x) dx  ,
$$
where $\chi_{f(x)\leq c_k(\I)}(x)$ denotes the characteristic function of 
the set $f(x)\leq c_k(\I)$. \\
Since $2(c(\I)+\varepsilon -f(x))+{\log k \over k} \geq \varepsilon$, the integrand is
dominated by a constant on $[0,2\pi]$. Therefore, by the Dominated
Convergence Theorem, we obtain
$$
\lim_{k \rightarrow +\infty} 
\int_0^{2\pi} 
{1 \over \sqrt{2(c(\I)+\varepsilon -f(x))+{\log k \over k}}}\chi_{f(x)\leq c_k(\I)}dx =
\int_0^{2\pi}{1 \over \sqrt{2(c(\I)+\varepsilon -f(x))}}dx =T(c(\I)+\varepsilon)  ,
$$
where $T(c(\I)+\varepsilon)$ is the period of
$$
\dot x=\gamma_0(c(\I)+\varepsilon,x) .
$$
Hence,  for any $\epsilon>0$, there 
exists $\tilde K_3(\I,\epsilon)$ such that, for any $k>\tilde K_3(\I,\epsilon)$ 
we have
$$
\int_0^{2\pi} 
{1 \over \sqrt{2(c(\I)+\varepsilon -f(x))+{\log k \over k}}}\chi_{f(x)\leq
  c_k(\I)}dx \geq {1\over 2} T(c(\I)+\epsilon) ,
$$
and therefore, for any $\epsilon>0$ if $k>\max \{ \tilde K_1(a_1,a_2),
\tilde K_2(\I,\epsilon), \tilde K_3(\I,\epsilon)\}$, 
$$
T_k \geq {a_1\over 2}\ T(c(\I)+\varepsilon) .
$$
Since $c(\I)=\max f$, one has
$$
\lim_{\epsilon \rightarrow 0^+}T(c(\I)+\varepsilon) =+\infty.
$$
Therefore, for any $\eta>0$ there exists $\epsilon(\eta)$ such that
$$
{a_1\over 2}\ T(c(\I)+\varepsilon(\eta)) > \eta ,
$$
and also for any $k>\max \{ \tilde K_1(a_1,a_2), 
\tilde K_2(\I,\epsilon(\eta)), \tilde K_3(\I,\epsilon(\eta))\}$, we have
$$
T_k \geq {a_1\over 2}\ T(c(\I)+\varepsilon(\eta))> \eta   .
$$
Hence we have proved
$$
\lim_{k\rightarrow +\infty} T_k=+\infty  .
$$
\hfill $\Box$ \\ \\
The rest of the proof will be formulated for $f(\varphi )=-\cos(\varphi)$.    
\begin{Lem}\label{lemlim2}
Let $f(\varphi)=-\cos \varphi$ and $\I \in (0,+\infty)$ be such that $c(\I) = \max f = 1$. Then, for 
any $t\in {\Bbb R}$, we have
$$
\lim_{k\rightarrow +\infty} x_k(t)=x(t) .
$$
\end{Lem}

\noindent
\textit{Proof.} Let us denote $d_k(t)= \norm{x_k(t)-x(t)}$. In order to overcome the lack of differentiability of $d_k$, for a constant $r > 0$, we introduce $e_k(t)= \sqrt{(x_k(t)-x(t))^2 + r^2}$, whose time derivative $\dot e_k(t)$ satisfies
\begin{eqnarray*}
\dot e_k(t) &\leq& \norm{\gamma_k(c_k(\I),x_k) - \gamma_0(c(\I),x)} \leq \\
&\leq& \norm{\gamma_k(c_k(\I),x_k) - \gamma_0(c_k(\I),x_k)} +
\norm{\gamma_0(c_k(\I),x_k)-\gamma_0(c(\I),x_k)} + 
\norm{\gamma_0(c(\I),x_k)-\gamma_0(c(\I),x)}  .
\end{eqnarray*}
In the sequel we denote $c_k:= c_k(\I),\, c:=c(\I)$ and we fix $t>0$.  By the uniform 
convergence of $\gamma_k$ to
$\gamma_0$, for any $\varepsilon>0$ there exists 
$\tilde K_1(\epsilon)$ such that for any $k> \tilde K_1(\epsilon)$, we have
$$
\norm{\gamma_k(c_k,x_k(\tau))-\gamma_0(c_k,x_k(\tau))} \leq \varepsilon 
$$
for any $\tau \in {\Bbb R}$, and specifically for any $\tau \in [0,t]$. \\
From Lemma \ref{lemlim1} and the convergence of $c_k$ to $c$, there also exist $\tilde K_2(t)$, $\rho(t)>0$
 such that for $k>{\tilde K_2}(t)$ one has: $T_k/2 > t$, 
 $ \norm{x_k(t)-\pi} > \rho(t)$ and also $f(x_k(\tau))<c_k$ for any $\tau \leq t$.
 Therefore, there exists $\lambda(t)<\infty$ such that
$$
\sup_{k}\sup_{\tau \leq t} {1\over \gamma_0(c,x_k(\tau))}=\lambda(t).
$$
As a consequence, we have
$$
\norm{\gamma_0(c_k,x_k(\tau))-\gamma_0(c,x_k(\tau))}={\norm{2(c_k-c)}\over 
\gamma_0(c_k,x_k(\tau))+\gamma_0(c,x_k(\tau))}\leq {\norm{2(c_k-c)}\over 
\gamma_0(c,x_k(\tau))} \leq 
 2 \lambda(t)\norm{c_k-c} .
$$
By the convergence of $c_k$ to $c$, for for any $\epsilon>0$ 
there exists ${\tilde K}_3(\epsilon)$ such that, for any 
$k> {\tilde K}_3(\varepsilon)$, $\norm{2(c_k-c)} \leq \epsilon$, and therefore for any $k>\max \{
{\tilde K_2}(t),{\tilde K}_3(\varepsilon)\}$ and any $\tau \in [0,t]$, 
$$
\norm{\gamma_0(c_k,x_k(\tau))-\gamma_0(c,x_k(\tau))} 
\leq \lambda(t) \varepsilon  .
$$
Moreover, since $1$ is a Lipschitz constant for $f$, we have 
$$
\norm{\gamma_0(c,x_k(\tau))-\gamma_0(c,x(\tau))} \leq {\norm{2(f(x_k(\tau))-f(x(\tau)))}
\over \gamma_0(c,x_k(\tau))+\gamma_0(c,x(\tau))}\leq {\norm{2(x_k(\tau)-x(\tau))}\over 
\gamma_0(c,x_k(\tau))}
\leq 2 \lambda(t) d_k(\tau) .
$$
Therefore, for any $\tau\leq t$, for any $\varepsilon>0$ and any 
$k>\max \{ {\tilde K}_1(\epsilon ),{\tilde K}_2(t),
{\tilde K}_3(\varepsilon)\}$, one has
$$
\dot e_k(\tau)\leq (1+\lambda(t))\varepsilon + 2 \lambda(t) d_k(\tau) .
$$
Since $d_k<e_k$, by the Gronwall Lemma, we have
$$
e_k(\tau)\leq {(1+\lambda(t))\varepsilon\over 2\lambda(t)}
\Big (e^{2\lambda(t)\tau}-1\Big ) ,
$$
and also 
$$
d_k(\tau)\leq {(1+\lambda(t))\varepsilon\over 2\lambda(t)}
\Big (e^{2\lambda(t)t}-1\Big ).
$$
For any $\eta>0$, let us consider $\epsilon(\eta)$ such that
${(1+\lambda(t))\varepsilon\over 2\lambda(t)}
\Big (e^{2\lambda(t)t}-1\Big )=\eta$. Then, for any $\eta>0$ 
and any $k>\max \{ {\tilde K}_1(\epsilon(\eta)),{\tilde K}_2(t),
{\tilde K}_3(\varepsilon(\eta))\}$, one has
$$
d_k(\tau) \leq e_k(\tau) \le \eta,
$$
that is
$$
\lim_{k\rightarrow \infty} d_k(\tau)=0 .
$$
Therefore, the function $d_k(t)$ converges pointwise to $0$ for any 
$t\in {\Bbb R}$, and the lemma is proved. \hfill $\Box$

 \begin{Lem}\label{lemlim3}
Let $f(x)=-\cos x$, and $\I$ such that $c(\I)=1$. Then, for 
any $\tilde t\in (0,1/2)$, we have
$$
\lim_{k\rightarrow +\infty} x_k(\tilde t\ T_k)=\pi .
$$
\end{Lem}

\noindent
\textit{Proof.} Since
$$
\lim_{t\rightarrow +\infty} x(t)=\pi
$$
for any $\varepsilon >0$ there exists $T(\varepsilon)$ such for $t\geq
T(\varepsilon)$ one has
$$
\norm{x(t)-\pi}\leq \varepsilon .
$$
Let us now consider $t=T(\varepsilon)$. By Lemma \ref{lemlim2} we have
$$
\lim_{k\rightarrow +\infty} x_k(t)=x(t) . 
$$
For any $\varepsilon>0$ there exists ${\tilde K}_1(\varepsilon)$ such that, 
for any $k>{\tilde K}_1(\varepsilon)$, we have
$$
\norm{x_k(T(\epsilon))-\pi}\leq 2 \varepsilon .
$$
Let us now fix $\tilde t\in (0,1/2)$, and consider 
${\tilde K}_2(\tilde t, \varepsilon)$ such that for any 
$k> K_2(\tilde t,\varepsilon)$, $\tilde{t} T_k>t$ 
(this is possible since by Lemma \ref{lemlim1} $T_k$ is a divergent 
sequence). Since $x_k(\tau)$ is monotone, from the inequalities $t < \tilde t T_k < T_k/2$, 
we obtain
$$
x_k(t) \leq x_k(\tilde tT_k) \leq x_k(T_k/2)=\pi ,
$$
and therefore, for any $k> \max \{ {\tilde K}_1(\varepsilon), 
K_2(\tilde t,\varepsilon)\}$,
$$
\norm{x_k(\tilde t T_k)-\pi}\leq 2 \varepsilon .
$$
The lemma is therefore proved. \hfill $\Box$
\vskip 0.4 cm
\noindent We can now prove the Proposition \ref{sigma-lim}. We consider the 
integral
$$
\int_{-\pi}^{\pi} u(x)\sigma_k(x)dx .
$$
and we change the integration variable from $x$ to $t$ 
by using $x=x_k(t)$, and then from $t$ to $\tilde t=t/T_k$, thus obtaining
$$
\int_{-\pi}^{\pi} u(x)\sigma_k(x)dx=
{1\over T_k}\int_{-{T_k\over 2}}^{T_k\over 2} u(x_k(t))dt 
=\int_{-{1\over 2}}^{1\over 2} u(x_k(\tilde t\ T_k))d\tilde t .
$$
We observe that, for any $\tilde t \in (0,1/2)$,
$$
\lim_{k\rightarrow +\infty}  u(x_k(\tilde t\ T_k))=u(\pi) .
$$
Also, for any $t\in (-1/2,0)$,
$$
\lim_{k\rightarrow \infty} u(x_k(\tilde t\ T_k))= u(-\pi)=u(\pi) .
$$
Therefore, the Lemma follows by the Dominated Convergence Theorem. \hfill $\Box$

\section{The Lambert function $W$} \label{W}
The Lambert function $W$ is defined as the multivalued function defined implicitly by the relation: 
\begin{equation}
z = W(z) e^{W(z)} 
\label{lambertfc}
\end{equation} 
for any complex number $z$. We only consider $W$ for $z\in[0,+\infty)$, so that it becomes single valued. In particular, $W(z) \geq 0$ for $z\in[0,+\infty)$ and $W$ is an increasing function. \\
\indent The asymptotic properties of $W$ may be characterized by asymptotic developments. We refer to \cite{cor1} and \cite{cor2} for all details and proofs. From these developments, one immediately obtains
\begin{equation}
\lim_{z\rightarrow +\infty} 
{W(z)\over \log z -\log \log z}=1
\label{liminf}
\end{equation}
and also
\begin{equation}
\lim_{z\rightarrow 0^+} 
{W(z)\over z-z^2}=1  .
\label{lim0}
\end{equation}
Formulas (\ref{lambertfc}), (\ref{liminf}) and (\ref{lim0}) are used extensively 
throughout this paper to define and prove the asymptotic properties of the key functions introduced in Definition \ref{def fond}. 

\vskip 0.4 cm
\noindent
{\bf Acknowledgments.}  We thank Francis Sullivan who helped us in improving our manuscript with many  useful suggestions. This research has been supported by project 
CPDA092941/09 of the University of Padova.

\end{document}